\documentclass[lettersize,journal]{IEEEtran}
\usepackage{amsmath,amsfonts}
\usepackage{algorithmic}
\usepackage{algorithm}
\usepackage{array}
\usepackage[caption=false,font=normalsize,labelfont=sf,textfont=sf]{subfig}
\usepackage{textcomp}
\usepackage{stfloats}
\usepackage{url}
\usepackage{verbatim}
\usepackage{graphicx}
\usepackage{cite}
\usepackage{derivative}
\usepackage{caption2}
\usepackage{bm}
\usepackage{amssymb}
\usepackage{xcolor}
\usepackage{amsthm}
\usepackage[pagewise]{lineno}

\usepackage{comment}

\newenvironment{list4}{
	\begin{list}{$\bullet$}{%
			\setlength{\itemsep}{0.05cm}
			\setlength{\labelsep}{0.2cm}
			\setlength{\labelwidth}{0.3cm}
			\setlength{\parsep}{0in} 
			\setlength{\parskip}{0in}
			\setlength{\topsep}{0in} 
			\setlength{\partopsep}{0in}
			\setlength{\leftmargin}{0.16in}}}
	{\end{list}}

\newtheorem{lemma}{Lemma}
\newtheorem{theorem}{Theorem}
\newtheorem{definition}{Definition}

\begin{document}
\switchlinenumbers  

\title{Estimating Sequences with Memory for Minimizing Convex Non-smooth Composite Functions}
\author{Endrit Dosti, 
Sergiy A. Vorobyov, and Themistoklis Charalambous

\thanks{E. Dosti and S. A. Vorobyov are with the Department of Information and Communications Engineering, Aalto University, Espoo, Finland (e-mails: endrit.dosti@aalto.fi and svor@ieee.org). Themistoklis Charalambous is with the Department of Electrical and Computer Engineering, University of Cyprus, Nicosia, Cyprus (e-mail: charalambous.themistoklis@ucy.ac.cy) and a Visiting Professor with the Department of Electrical Engineering and Automation, Aalto University, Espoo, Finland. {An earlier version of this paper was presented at the IEEE Conference on Decision and Control, Cancun, Mexico, Dec. 2022 [DOI:10.1109/CDC51059.2022.9993313]}}

}


\maketitle
\begin{abstract}
First-order optimization methods are crucial for solving large-scale data processing problems, particularly those involving convex non-smooth composite objectives. For such problems with convex non-smooth composite objectives, we introduce a new class of generalized composite estimating sequences, devised by exploiting the information embedded in the iterates generated during the minimization process. Building on these sequences, we propose a novel accelerated first-order method tailored for such objective structures. This method features a backtracking line-search strategy and achieves an accelerated convergence rate, regardless of whether the true Lipschitz constant is known. Additionally, it exhibits robustness to imperfect knowledge of the strong convexity parameter, a property of significant practical importance.
The method's efficiency and robustness are substantiated by comprehensive numerical evaluations on both synthetic and real-world datasets, demonstrating its effectiveness in data processing applications.
\end{abstract}


\section{Introduction}
\label{sec: Introduction}
\IEEEPARstart{R}{ecent} research in first-order methods for solving large-scale data processing problems has been largely focused on exploring different approaches to the acceleration of gradient-based methods \cite{Pesquet}. For the problem of minimizing smooth convex functions\footnote{Minimization of smooth convex functions is a problem of theoretical interest although its applicability in signal processing is quite limited compared to minimization of non-smooth convex functions, which we consider here.}, we recently developed a method by extending the estimating sequences framework \cite{Endrit1} that
converges faster than the Fast Gradient Method (FGM) \cite{Nesterov_83, Nesterov_book}. In yet another framework, the continuous-time limit of FGM has been modeled as a second-order differential equation \cite{pmlr-v40-Flammarion15, Su_Boyd_Candes, WibisonoE7351}. In another newly developed framework \cite{drori-teboulle}, the authors have cast the improvement of the worst-case behavior of an algorithm as an optimization problem. Based on this framework, an optimal method for minimizing smooth convex functions has been presented in \cite{taylor2021optimal}. Despite the promising theoretical analysis, the applicability of these methods in the current form is restricted only to minimizing smooth convex functions, and their generalization capabilities remain unclear.

Considering the different strategies that have been developed for accelerating gradient-based methods, estimating sequence methods continue to play a central role in the field (see \cite{d2021acceleration} and references therein). First, for the case of differentiable convex functions, such methods are optimal in the sense of \cite{Nemirovski_Yudin}, {that is, such first-order methods are optimal (with accuracy to a multiplicative constant) in terms of the required number of iterations for achieving a given tolerance.} Second, they are efficient in practice and can work well with backtracking line-search \cite{Iulian_1, nesterov2015universal}. Third, they can be used to devise fast second-order and higher-order methods \cite{fast_newton_method, nesterov2020inexact}. Fourth, their efficiency has also been established in the context of applications to distributed optimization, nonconvex optimization, stochastic optimization, and many more (see \cite{jakovetic2014fast, ghadimi2016accelerated, kulunchakov:hal-01993531, ahn2020nesterov, li2020revisit} and the references therein). As discussed in \cite{Nesterov_book}, different estimating sequences can be used to enable the accumulation of global information of the objective function. One of the main challenges with the framework is the design of estimating functions that are used to construct the estimating sequences. 

{The estimating sequences framework has been formalized in \cite{Nesterov_book, baes2009estimate}.} For the broader class of minimizing convex functions with composite structure, which is important to this paper, a popular method is the Accelerated Multistep Gradient Scheme (AMGS) \cite{Nesterov_2007}, which exhibits an accelerated convergence rate. The method has the disadvantage of requiring two projection-like operations per iteration, which translates into an increased runtime of the method and inhibits its deployment to practical large-scale optimization setups \cite{nesterov2014subgradient}. Another popular method is the Fast Iterative Shrinkage-Thresholding Algorithm (FISTA) \cite{FISTA}. Unlike AMGS, it requires one projection-like operation per iteration and has been proven to exhibit an accelerated convergence rate. Nevertheless, as we will also see in the numerical section, the method converges slower than AMGS. At first glance, FISTA does not appear to be an estimating sequence method. Nevertheless, links between FISTA and estimating sequence methods have been established in \cite{Iulian_2}. {In \cite{Dosti_1}, the authors have introduced COMET, which is a new first-order accelerated algorithms built based on the estimating sequences framework used for devising FGM. Similar to FISTA, the method proposed therein requires one projection-like operation per iteration, it is more efficient than AMGS, and leads to a theoretically established constant improvement of the convergence rate.}


Contrasting the analysis conducted for AMGS in \cite{Nesterov_2007} with FGM in \cite{Nesterov_book}, we can see that different estimating functions were used. As discussed earlier, the lack of uniqueness of the estimating sequences is one of the main challenges in developing methods under such a framework. In theory, 
all the aforementioned methods exhibit an accelerated convergence rate, but they perform differently in practice. {Moreover, in \cite{Endrit1, Dosti_2}, the authors have shown how to devise generalized estimating sequences, which can be used to construct faster algorithms. The generalized estimating sequences accommodate additional terms, which represent any additional knowledge about the objective. In the black-box framework, the only additional knowledge about the objective available is the memory of the iterative updates, which also allows to better estimate the curvature of the objective -- the information needed to accelerate the convergence and reduce algorithm's sensitivity to accurate knowledge of hyperparameters, such as Lipschitz constant and strong convexity parameter. The development in \cite{Endrit1, Dosti_2} is, however, limited to considering the minimization of smooth convex functions only.}

{As demonstrated above, the interest in first-order optimization algorithms is very significant, especially towards addressing the associated practical issues such as better convergence for real-world data processing, robustness to a number of hyperparameter of an algorithm such as Lipschitz constant and strong convexity parameter, relaxation of assumptions and improvement of convergence guarantees. New methods for practical performance improvement and software for performance estimation \cite{PEPit} are of high interest. The non-smooth objective structure is especially of high importance in applications to data and signal processing, where the problems that manifest themselves as minimization of composite non-smooth convex functions (rather than smooth convex functions as in \cite{Endrit1, Dosti_2}) are a lot more common. The composite objective structure with a non-smooth part appears for example in the problems with sparsity constraints, which are typically addressed by adding a non-smooth penalty to the objective.} 

{Focusing on improving convergence both theoretically and for real-world data processing, robustness to algorithm's hyperparameters, relaxation of assumptions and improvements of convergence guarantees, we propose here a generalization of estimating sequences for composite non-smooth objectives and develop new algorithm and new convergence results for the algorithm.\footnote{{Preliminary results towards such an extension of generalized estimating sequences framework to non-smooth composite objectives were reported in \cite{Dosti_3}.}} This work can be viewed as a significant advancement of the framework, which we first introduced in \cite{Endrit1}, to the practically more appealing case of non-differentiable objectives. This case is technically harder but significantly widens the applicability for addressing signal processing problems. }
%
The main contributions of the article are as follows.
\begin{list4}
	\item We introduce a new structure for the estimating functions, which we call the \emph{generalized composite estimating functions}. The proposed estimating functions are devised by making use of the following: \textit{(i)}~A new term created by adding the previously constructed estimating functions; \textit{(ii)}~The gradient mapping framework \cite{Nemirovski_Yudin}; \textit{(iii)}~A tighter lower bound on the objective function. 
	\item Using our proposed estimating sequences, we devise a new accelerated method for minimizing convex non-smooth composite functions. Moreover, we present an efficient line-search strategy which is used to estimate the step size. Our proposed method requires only one projection-like operation per iteration, which is lower than the respective requirement for AMGS. 
	\item We prove that our proposed method exhibits an accelerated convergence rate despite the inaccurate knowledge of the Lipschitz constant. Note that in practice, it is reasonable to assume that the computational cost of finding an upper bound to the Lipschitz constant is acceptably low, but it is unreasonable to assume exact knowledge of it. 
	\item We prove that the way our proposed method is initialized is robust to the inaccurate knowledge of the strong convexity parameter, 
	which further reduces the additional computational burden of having to estimate such a parameter. Indeed, there exists no cost-efficient generic approach for estimating the strong convexity parameter. 
	\item We demonstrate the efficiency of our proposed method as compared to the existing benchmarks. Using real-world datasets, in our computational experiments, we also highlight the robustness of our proposed method in cases when the strong convexity parameter and Lipschitz constant 
	are not known, which is important in practical data processing applications where these parameters are unknown and very computationally expensive to estimate. 
\end{list4}

The article is organized as follows. Section~\ref{prel} defines the setup and the necessary preliminaries. In Section~\ref{sec: Proposed Method}, we present the generalized composite estimating sequences and show how they can be used to build our proposed method. In Section~\ref{sec: Convergence Analysis}, we prove the convergence results for our proposed method. In Section~\ref{sec: Simulations}, we depict the numerical performance of our proposed method and compare it with several existing benchmarks. We consider several types of optimization problems and demonstrate the efficiency of our proposed method. Last, in Section~\ref{Con}, we summarize the main findings of the paper. All technical proofs are given in the Appendix.

\section{Preliminaries}
\label{prel}
In the sequel, we will focus on devising an accelerated black-box method for solving convex optimization problems with composite objective functions. The typical structure for such problems is
\begin{equation}
	\begin{aligned}
		\label {opt_prob}
		F({\boldsymbol x}) = f({\boldsymbol x}) + \tau g({\boldsymbol x}), \quad \tau > 0,
	\end{aligned}
\end{equation}
where $f: \mathcal{R}^n \rightarrow \mathcal{R}$ is a differentiable convex function and $g: \mathcal{R}^n \rightarrow \mathcal{R}$ is a simple convex lower semi-continuous function. Here, $\mathcal{R}$ and $ \mathcal{R}^n $ are the set of real numbers and the set of real-valued vectors of size $n \times 1$, respectively. {Such composite structure often appears in signal and image processing problems when the regularizer is not a smooth function. For example in LASSO, the objective is the least squares and the regularizer is the $l_1$-norm of the vector of optimization variables.} The simplicity of $g$ implies that the complexity of computing the proximal map
\begin{equation}
	\begin{aligned}
		\label {prox_g}
		\text{prox}_{\tau g} \triangleq \; & \text{arg} \underset{{\boldsymbol z} \in {\mathcal{R}^n}}{\min}
		& &\left(g({\boldsymbol z}) + \frac{1}{2 \tau} \|{\boldsymbol z} - {\boldsymbol x} \|^2\right), \quad {\boldsymbol x} \in \mathcal{R}^n, 
	\end{aligned}
\end{equation}
is $\mathcal{O} (n)$ \cite{Boyd_prox_alg}. Herein $\| \cdot \|$ denotes the $l_2$ norm of a vector.

Assuming that $g({\boldsymbol x})$ has strong convexity parameter $\mu_g \geq 0$, we use the following strong convexity transfer
\begin{align} \label{FFF}
	F({\boldsymbol x}) \! &= \! \underbrace{ \left( \! f({\boldsymbol x}) + \frac{\tau \mu_g}{2} \| {\boldsymbol x} - {\boldsymbol x}_0 \|^2 \! \right) }_{\hat{f} ({\boldsymbol x})} \! + \tau \! \underbrace{ \left( \! g({\boldsymbol x}) - \frac{\mu_g}{2} \| {\boldsymbol x} - {\boldsymbol x}_0 \|^2 \! \right) }_{\hat{g} ({\boldsymbol x})}
\end{align}
to facilitate the tractability of the derivations presented in the sequel. Here ${\boldsymbol x}_0$ is an initial value of ${\boldsymbol x}$. Based on \eqref{FFF}, we have $L_{\hat{f}} = L_f + \tau \mu_g$, $\mu_{\hat{f}} = \mu_f + \tau \mu_g$ and $\mu_{\hat{g}} = 0$. Here, $L_f$ and $\mu_f$ are the Lipschitz constant and the strong convexity parameter of $f$, respectively.

For all ${\boldsymbol y}, 
\, {\boldsymbol x} \in \mathcal{R}^n$, and $L \geq L_{\hat{f}}$, let us define
\begin{align}
	m_{L} ({\boldsymbol y}; {\boldsymbol x}) \! \triangleq \! \hat{f}({\boldsymbol y}) \! + \! \nabla \hat{f}^T ({\boldsymbol y}) ({\boldsymbol x} \! - \! {\boldsymbol y}) \! + \! \frac{L}{2} \| {\boldsymbol x} \!-\! {\boldsymbol y} \|^2 \! + \! \tau \hat{g}({\boldsymbol x}), \label{m_l}
\end{align}
where $\nabla$ denotes gradient and $(\cdot)^T$ stands for transposition. The following bounds for $\hat{f}({\boldsymbol x})$ and $\hat{g}({\boldsymbol x})$ will be useful in the analysis
\begin{align}
	\hat{f}({\boldsymbol x}) &\leq \hat{f}({\boldsymbol y}) + \nabla \hat{f}^T ({\boldsymbol y}) ({\boldsymbol x} - {\boldsymbol y}) + \frac{L_{\hat{f}}}{2} \| {\boldsymbol y} - {\boldsymbol x} \|^2, \label{upper_bound} \\
	\hat{g}({\boldsymbol x}) &\geq \hat{g}({\boldsymbol y}) + s^T ({\boldsymbol y}) ({\boldsymbol x} - {\boldsymbol y}).  \label{lower_bound_g}
\end{align}
Then, considering \eqref{m_l} and \eqref{upper_bound}, we have
\begin{align}
	\label{LLL_bound}
	m_{L}({\boldsymbol y}; {\boldsymbol x}) \geq F({\boldsymbol x}), \quad \forall {\boldsymbol x}, {\boldsymbol y} \in \mathcal{R}^n.
\end{align} 

Next, we define the composite gradient mapping as {\cite{Nesterov_book}}
\begin{align}
	\label{T_L(y)}
	T_{L}({\boldsymbol y}) &\triangleq \arg \underset{{\boldsymbol x} \in \mathcal{R}^n}{\min} \; m_{L} ({\boldsymbol y}; {\boldsymbol x}).
\end{align}
Then, the reduced composite gradient is defined as
\begin{align}
	\label{r_l}
	r_{L} ({\boldsymbol y}) \triangleq L \left({\boldsymbol y} - T_{L}({\boldsymbol y})\right).
\end{align}

Consider now the optimality conditions for \eqref{T_L(y)} {\cite{Nesterov_book}}:
\begin{align}
	\nonumber
	\nabla m_L^T ({\boldsymbol y};T_L({\boldsymbol y})) ({\boldsymbol x} - T_L({\boldsymbol y})) &\geq 0, \\ \label{r(y)}
	\big(\nabla f({\boldsymbol y}) + L(T_L ({\boldsymbol y}) - {\boldsymbol y}) + \tau s_L({\boldsymbol y}) \big)^T({\boldsymbol x} - T_L({\boldsymbol y})) &\geq 0,
\end{align}
where $s_L ({\boldsymbol y}) \in \partial g(T_L({\boldsymbol y}))$ is a subgradient and $\partial g(T_L({\boldsymbol y}))$ is the subdifferential. In \eqref{r(y)}, let
\begin{align}
	\label{iu}
	\nabla f({\boldsymbol y}) + L(T_L ({\boldsymbol y}) - {\boldsymbol y}) + \tau s_L({\boldsymbol y}) = 0.
\end{align}
Substituting \eqref{iu} in \eqref{r_l} yields
\begin{align}
	\label{reduced_grad}
	r_L({\boldsymbol y}) = L({\boldsymbol y} - T_L({\boldsymbol y})) = \nabla f({\boldsymbol y}) + \tau s_L({\boldsymbol y}).
\end{align}

In multiple places in the paper, we will make use of the following tight lower bound on the objective function, which we first established in \cite[Theorem~1]{Dosti_1, Dosti_2}.
\begin{theorem}
	\label{thm2}
	Let $F({\boldsymbol x})$ be a composition of an $L_{\hat{f}}$-smooth and $\mu_{\hat{f}}$-strongly convex function $\hat{f}({\boldsymbol x})$, and a simple convex function $\hat{g}({\boldsymbol x})$, as given in \eqref{FFF}. For $L \geq L_{\hat{f}}$, and ${\boldsymbol x}, \, {\boldsymbol y} \in \mathcal{R}^n$ we have 
	\begin{align}
		\label{13}
		F({\boldsymbol x}) &\geq \hat{f}(T_L({\boldsymbol y})) + \tau \hat{g} (T_L({\boldsymbol y})) + r_L^T ({\boldsymbol y}) \left({\boldsymbol x} - {\boldsymbol y} \right) \nonumber \\
		&+ \frac{\mu_{\hat{f}}}{2} \| {\boldsymbol x} - {\boldsymbol y} \|^2 + \frac{1}{2L} \|r_L({\boldsymbol y})\|^2 .
	\end{align}
\end{theorem} 

For establishing the bound on the convergence rate or a required number of iterations for achieving a given tolerance, we will also make use of the following upper bound on the difference $F({\boldsymbol x}_0) - F({\boldsymbol x}^*)$, where $F({\boldsymbol x}_0)$ and $F({\boldsymbol x}^*)$ are the objective function values at the starting point ${\boldsymbol x}_0$ and at optimality ${\boldsymbol x}^*$, respectively. 
\begin{theorem}	
	Let $F({\boldsymbol x})$ be a convex function with composite structure as shown in \eqref{opt_prob}. 
	Then, for any feasible starting point ${\boldsymbol x}_0$, we have
	\begin{align}
		\label{ll} 
		F({\boldsymbol x}_0) - F({\boldsymbol x}^*) & \leq \frac{L_0}{2} \|{\boldsymbol x}_0 - {\boldsymbol x}^*\|^2, 
	\end{align}
	where {$L_0$ is the estimate of the value of $L$ at iteration $k=0$}, that is, the Lipschitz constant of $F({\boldsymbol x})$ at the starting point ${\boldsymbol x}_0$ or its upper bound. 
\end{theorem}
The inequality \eqref{ll} is straightforward for smooth functions, but it requires a tedious proof for convex functions with composite structure, and it can be found in \cite[Lemma~5]{Dosti_2}.

Finally, note that the goal of a numerical optimization scheme is to devise a sequence of iterates ${\boldsymbol x}_0, {\boldsymbol x}_1, \ldots {\boldsymbol x}_k$, which goes arbitrarily close to the optimal solution ${\boldsymbol x}^*$ (within some tolerance $\epsilon > 0$). The set $\mathcal{R}^n$ is much larger than the search area of interest at iteration $k$, given as $\left\{{\boldsymbol x} \,|\, {\boldsymbol x}_0 + \text{span} \{\nabla f({\boldsymbol x}_0), \ldots \nabla f ({\boldsymbol x}_{k-1})\} \right\}$ for designing any first-order method. Here $\text{span} \{ \cdot \}$ denotes a span of a set of vectors. Thus, instead of considering the largest possible set in which the objective function is defined, i.e., $\mathcal{R}^n$ for unconstrained optimization, we will establish results in our paper for the 
subset of $\mathcal{R}^n$ needed for the method that we will design, i.e., 
\begin{align}
	\label{defQ} 
\mathcal{Q} = \left\{{\boldsymbol x} \,|\, {\boldsymbol x}_0 + \text{span} \{\nabla f({\boldsymbol x}_0), \ldots \nabla f ({\boldsymbol x}_{k-1}), \ldots\} \right\} \subset \mathcal{R}^n.
\end{align}

\section{Proposed Method}
\label{sec: Proposed Method}
Consider the following definition for the generalized composite estimating sequences. 
\begin{definition}
	\label{def__1}
	The sequences $\{\Phi_{k}\}_{k}$ and $\{\lambda_{k}\}_{k}$, $\lambda_{k} \geq 0$, are called generalized composite estimating sequences of the function $F(\cdot)$ defined in \eqref{FFF}, if there exists a sequence of bounded functions $\{\psi_k\}_k$, $\lambda_{k} \rightarrow 0$ as $k \rightarrow \infty$, and $\forall {\boldsymbol x} \in \mathcal{Q}$, $\forall k = 0, 1, \cdots$ we have
	\begin{equation}
		\label{def_1}
		\Phi_{k} ({\boldsymbol x}) \leq\lambda_{k} \Phi_{0}({\boldsymbol x}) + (1 - \lambda_{k}) \left(F({\boldsymbol x}) - \psi_k ({\boldsymbol x})\right).
	\end{equation}
\end{definition} 
{Note that in \eqref{def_1}, we have an additional term $\psi_k ({\boldsymbol x})$ as compared to the definition of standard composite estimating sequences \cite{Nesterov_book, Dosti_1} that, if chosen carefully, can impact the convergence of the iterates of the corresponding optimization algorithm. Thus, our objective in the sequel is to demonstrate a concrete design for this term and further demonstrate both analytically and numerically the improvement in convergence rate of the corresponding optimization algorithm. We aim to stay within the black-box setup, that is, we have no prior knowledge about the particular structure of the objective, except that it is a composite function as defined by \eqref{opt_prob}, and thus, the only additional information for constructing $\psi_k ({\boldsymbol x})$ is the history of iterative updates.}

Let us now use the above defined generalized composite estimating sequences to characterize the convergence rate of the minimization process summarized in the form of the following lemma.
\begin{lemma}
	\label{SFGM_lemma_1}
	If for a sequence $\{{\boldsymbol x}_k\}_{k}$ we have $F({\boldsymbol x}_k) \leq \Phi_{k}^*  \triangleq \underset{{\boldsymbol x} \in {\mathcal{Q}}}{ \min } \Phi_{k} ({\boldsymbol x})$, then 
    \begin{equation}
        F({\boldsymbol x}_k) - F({\boldsymbol x}^*) \leq \lambda_{k} \left[ \Phi_{0}({\boldsymbol x}^*) - F({\boldsymbol x}^*) \right]- (1 - \lambda_k) \psi_k ({\boldsymbol x}^*) ,
    \end{equation}
        where ${\boldsymbol x}^* \triangleq \arg \underset{{\boldsymbol x} \in \mathcal{Q}}{\min} F({\boldsymbol x})$.
\end{lemma}
For the proof see Appendix {\it Proof of Lemma~\ref{SFGM_lemma_1}.}

Let us now present the estimating functions that will be used to devise our proposed method. 
\begin{lemma}
	\label{SFGM_lemma_2}
	Assume that there exist sequence $\{\alpha_k\}_{k}$, where $\alpha_{k} \in (0, 1)$ $\forall k$, such that $\sum_{k = 0}^{\infty} \alpha_{k} = \infty$; sequence $\{\psi_k\}_k$ with an upper bound $\Psi$, such that $\{\psi_k\}_k \geq 0$; and an arbitrary sequence $\{{\boldsymbol y}_k\}_{k}$. Furthermore, let $\psi_0 ({\boldsymbol x}) = 0$, $\lambda_{0}$ = 1 and assume that the estimates $L_k, \, \forall k$ of the Lipschitz constant $L_{\hat{f}}$ are selected in a way that inequality \eqref{upper_bound} is satisfied for all the iterates ${\boldsymbol x}_k$ and ${\boldsymbol y}_k$. Then, the sequences $\{\Phi_{k}\}_{k}$ and $\{\lambda_{k}\}_{k}$, which are defined recursively as
	\begin{align}
		\label{lambda_recursive}
		\lambda_{k+1} &= (1 - \alpha_k) \lambda_{k}, \\ 
		\Phi_{k+1} ({\boldsymbol x})  &= (1 - \alpha_k) \left(\Phi_{k} ({\boldsymbol x}) + \psi_{k} ({\boldsymbol x}) \right) - \psi_{k+1} ({\boldsymbol x}) - \Psi \nonumber \\  &+  \! \alpha_{k} \left( F\left(T_{L_k} ({\boldsymbol y}_k) \right) + \psi_{k} ({\boldsymbol x}) + \frac{1}{2L_k} \|r_{L_k} ({\boldsymbol y}_k)\|^2 \right) \nonumber \\ &+ \alpha_k \left( r_{L_k}^T ({\boldsymbol y}_k) ({\boldsymbol x} - {\boldsymbol y}_k)  + \frac{\mu_{\hat{f}}}{2} \| {\boldsymbol x} - {\boldsymbol y}_k \|^2\right), \label{phi_k+1_SFGM}
	\end{align}
	are {generalized} composite estimating sequences.
\end{lemma}
For the proof see Appendix {\it Proof of Lemma~\ref{SFGM_lemma_2}.}

Let us now compare between the different estimating sequence constructions that exist in the literature. First, observe that the estimating sequences used to construct FGM in \cite[Lemma~2.2.4]{Nesterov_book} are the instance of our proposed generalized composite estimating sequences obtained when $\tau = 0$ and $\{\psi_k\}_k = 0$. Moreover, both types of estimating sequences can be used to measure the convergence rate of the minimization process. In this sense, the framework presented herein, is a generalization of the estimating sequences framework. Comparing our generalized composite estimating sequences to \cite{Dosti_1, Dosti_2}, the introduction of the terms $\{\psi_k\}_k$ has an additional impact on the convergence rate of the minimization process. 

There are different ways to choose $\{\Phi_k\}_{k}$ and $\{\psi_k 
\}_{k}$. Let {$\phi_k^* 
\in \mathcal{R}$
is the minimal value that the estimating function can take for ${\boldsymbol x} \in \mathcal{Q}$, where $\mathcal{Q}$ is given by \eqref{defQ},} $\gamma_k \in \mathcal{R}^+$ ($\mathcal{R}^+$ is the set of real non-negative numbers), ${\boldsymbol v}_k \in \mathcal{Q}$, 
$\forall k = 0, 1, \ldots$ and define the terms $\{ \Phi_k \}_{k}$ as
\begin{align}
	\label{phi}
	\Phi_k({\boldsymbol x}) \triangleq \phi_k^* + \frac{\gamma_k}{2} \| {\boldsymbol x} - {\boldsymbol v}_k \|^2 - \psi_k ({\boldsymbol x}), \; k = 1, 2, \ldots .
\end{align}
{Note that we select above a parabolic structure for $\Phi_k ({\boldsymbol x})$, where ${\boldsymbol v}_k$ has then a meaning of the center of the parabola. Since our goal is
to construct a generalized version of an accelerated algorithm for minimizing a composite objective \eqref{opt_prob} with no additional prior knowledge about objective's particular structure (black-box
setup), a simple and quite generic approach to designing $\psi_k ({\boldsymbol x})$ is to let the terms in the sequence $\{\Phi_k\}_{k}$ ‘‘self-regulate’’ based on the memory of algorithm's iterative updates. Particularly, the terms of the sequence $\{\psi_k \}_{k}$ can be chosen to account for the accumulation of the terms in the sequence $\{\Phi_k \}_{k}$ as follows }
\begin{align}
	\label{psii}
	\psi_k ({\boldsymbol x}) \triangleq \sum_{j = 1}^{k-1} \! \beta_{j,k} \frac{\gamma_{j}}{2} \| {\boldsymbol x} - {\boldsymbol v}_{j} \|^2, \quad k = 1, 2, \ldots,
\end{align}
where $\beta_{j,k} \in [0, 1], \; j = 1, \ldots, k-1$. 

Considering the definition introduced above for $\Phi_k({\boldsymbol x})$ and $\psi_k({\boldsymbol x})$, it is of interest to assess the conditions for $\psi_k ({\boldsymbol x})$ that ensure the convexity of $\Phi_k ({\boldsymbol x})$. Since both functions are twice differentiable, assessing the second order condition for \eqref{phi}, we have $\sum_{j = 1}^{k-1} \beta_{j,k} \gamma_{j} \leq \gamma_{k}$. Moreover, we also restrict $\sum_{j = 1}^{k-1} \beta_{j,k} \gamma_{j} \leq \mu$. Combining these conditions, we reach
\begin{align}
	\label{psi_bound_}
	\sum_{j = 1}^{k-1} \beta_{j,k} \gamma_{j} \leq \min \left(\gamma_{k}, \mu\right).
\end{align}

We can find the minimal value of the estimating function introduced in \eqref{phi} as
\begin{align}
	\label{Phi^*}
	\Phi_k^* &= \underset{{\boldsymbol x} \in {\cal Q}}{ \min } \Phi_{k} ({\boldsymbol x}) \nonumber \\
	&= \phi_k^* + \frac{\gamma_k}{2} \| {\boldsymbol x}_{\Phi_k}^* - {\boldsymbol v}_k \|^2 - \sum_{j = 1}^{k-1} \! \frac{\beta_{j,k} \gamma_{j}}{2} \| {\boldsymbol x}_{\Phi_k}^* - {\boldsymbol v}_{j} \|^2,
\end{align}
where ${\boldsymbol x}_{\Phi_k}^* \triangleq \arg \min_{{\boldsymbol x} \in \mathcal{Q}} \Phi_k ({\boldsymbol x})$. The values of the parameters still need to be computed in a recurrent manner. The following Lemma captures these relations for the components of $\{ \Phi_k \}_{k}$ introduced in \eqref{phi}.
\begin{lemma}
	\label{SFGM_lemma_3}
	Assume that the coefficients $\beta_{i, k}$ are selected such that \eqref{psi_bound_} is satisfied and let $\phi_{0} ({\boldsymbol x}) = \phi_{0}^* + \frac{\gamma_{0}}{2} \| {\boldsymbol x} - {\boldsymbol v}_{0}||^2$, where $\gamma_0 \in \mathcal{R}^+$ and ${\boldsymbol v}_0 = {\boldsymbol x}_0$, for example. Then, the process defined in Lemma \ref{SFGM_lemma_2} preserves the canonical form of the function $\Phi_{k} ({\boldsymbol x})$ presented in \eqref{phi}, where the sequences $\{\gamma_{k}\}_{k}$, $\{{\boldsymbol v}_{k}\}_{k}$ and $\{\phi_{k}^*\}_{k}$ can be computed as 
	\begin{align}
		\label{gamma_expr}
		\gamma_{k+1} &= (1-\alpha_k)\gamma_{k} + \alpha_k \sigma_k, \\
		\label{v_value}
		{\boldsymbol v}_{k + 1}  &= \frac{1}{\gamma_{k + 1}} \! \bigg((1 \! - \! \alpha_k) \gamma_k {\boldsymbol v}_{k} \! + \! \alpha_k \Big( \mu_{\hat{f}} {\boldsymbol y}_k + \sum_{j = 1}^{k-1} \beta_{j, k} \gamma_j {\boldsymbol v}_j \nonumber \\
		& - r_{L_k} ({\boldsymbol y}_k) \Big) \bigg), \\
		\phi_{k+1}^*  &= (1-\alpha_k) \phi_k^* + \! \alpha_{k} \xi_k , \label{phi_{k+1}^*}
	\end{align}
	where 
    \begin{equation}
        \label{sigmaK}
        \sigma_k \triangleq \mu_{\hat{f}} + \sum_{j = 1}^{k-1} \beta_{j, k} \gamma_j
    \end{equation} 
    and $\xi_k$ is defined in \eqref{xi_k} at the bottom of the next page.
\begin{figure*}[!b]
\begin{align} \label{xi_k}
	\xi_k  &\triangleq F\left(T_{L_k} ({\boldsymbol y}_k) \right) + \frac{1}{2L_k} \|r_{L_k} ({\boldsymbol y}_k) \|^2 + \sum_{j = 1}^{k-1} \beta_{j, k} \gamma_j \| {\boldsymbol y}_k - {\boldsymbol v}_j \|^2 \! - \! \frac{L_k^2 \alpha_k}{2 \gamma_{k+1}} \| {\boldsymbol y}_k \! - \! T_{L_k} \!({\boldsymbol y}_k) \|^2 \!+\! \frac{\gamma_k (1 - \alpha_k) \sigma_k}{2 \gamma_{k+1}} \| {\boldsymbol y}_k \! - \! {\boldsymbol v}_k \|^2  \nonumber \\
	&+ \frac{(1 \! - \! \alpha_k) \gamma_k}{\alpha_k \gamma_{k+1}} \| {\boldsymbol x}_{\Phi_k}^* - {\boldsymbol v}_k \|^2 + \sum_{j = 1}^{k} \frac{\beta_{j, k+1} \gamma_j}{2 \alpha_k} \| {\boldsymbol x}_{\Phi_{k+1}}^* - {\boldsymbol v}_j \|^2 \!+\! \frac{\alpha_k (1-\alpha_k)}{\gamma_{k+1}} \sum_{j = 1}^{k-1} \beta_{j, k} \gamma_j({\boldsymbol v}_j - {\boldsymbol y}_k)^T r_{L_k}({\boldsymbol y}_k) \nonumber \\ 
	&+ \frac{\alpha_k^2}{\gamma_{k+1}} \sum_{j = 1}^{k-1} \beta_{j, k} \gamma_j \| {\boldsymbol v}_j - {\boldsymbol y}_k \| \; \| r_{L_k}({\boldsymbol y}_k) \| \!+\! \frac{\gamma_k (1-\alpha_k)}{\gamma_{k+1}} \Big( ({\boldsymbol v}_k - {\boldsymbol y}_k)^T r_{L_k}({\boldsymbol y}_k) + \sum_{j = 1}^{k-1} \beta_{j, k} \gamma_j \| {\boldsymbol y}_k - {\boldsymbol v}_j \| \; \| {\boldsymbol y}_k - {\boldsymbol v}_k \| \Big).
\end{align}
\end{figure*}	
\end{lemma} 
For the proof see Appendix {\it Proof of Lemma~\ref{SFGM_lemma_3}.}

Comparing the result obtained in Lemma \ref{SFGM_lemma_3} with that of \cite[Lemma 2.2.3]{Nesterov_book}, it can be seen that the recursive relations obtained for computing the elements of $\{{\boldsymbol v}_{k}\}_{k}$ and $\{\phi_{k}^*\}_{k}$ now reflect on the usage of a new lower bound on the function that is being minimized, and the reduced composite gradient. Note that the recurrent relations for computing $\{\gamma_{k}\}_{k}$, $\{{\boldsymbol v}_{k}\}_{k}$ and $\{\phi_{k}^*\}_{k}$ all reflect the presence of the added memory term that was used to construct them. Comparing the above obtained results \cite{Dosti_1, Dosti_2}, we can observe the additional terms coming from the newly introduced memory terms into the generalized composite estimating sequences.

To devise our proposed method, we will use an inductive argument. Assume that for a step $k$ we have
\begin{align}
	\Phi_k^* \!\stackrel{\eqref{Phi^*}}{=}\! \phi_k^* \!+\! \frac{\gamma_k}{2} \| {\boldsymbol x}_{\Phi_k}^* \!-\! {\boldsymbol v}_k \|^2 \!-\! \sum_{j = 1}^{k-1} \! \frac{\beta_{j,k} \gamma_{j}}{2} \| {\boldsymbol x}_{\Phi_k}^* \! - \! {\boldsymbol v}_{j} \|^2 \geq F({\boldsymbol x}_k). \label{Phi_k^*}
\end{align}
For the inductive argument to be complete, we need to establish that $\Phi_{k+1}^* \geq F({\boldsymbol x}_{k+1})$. Considering the assumption for iteration $k$ in \eqref{Phi_k^*},
\eqref{phi_{k+1}^*} yields
\begin{align}
	\phi_{k + 1}^* &\geq (1 - \alpha_k) F({\boldsymbol x}_k) + \alpha_{k} \xi_k. \label{36}
\end{align}
Using \eqref{13} in \eqref{36}, we reach 
	\begin{align}
		\phi_{k+1}^*  &\geq (1 - \alpha_k) \bigg( F(T_{L_k}({\boldsymbol y}_k)) + r_{L_k}^T ({\boldsymbol y}_k) \left({\boldsymbol x}_k - {\boldsymbol y}_k\right) \nonumber \\
		&+ \frac{\mu}{2} \| {\boldsymbol x}_k - {\boldsymbol y}_k \|^2 + \frac{1}{2L_k} \| r_{L_k} ({\boldsymbol y}_k)||^2 \bigg) + \alpha_{k} \xi_k. \label{37}
	\end{align}
Substituting \eqref{xi_k} into \eqref{37} and performing some straightforward linear transformations, we get inequality \eqref{38} shown at the top of the next page.
\begin{figure*}
	\setcounter{equation}{31}
	\begin{align}
		\phi_{k+1}^*  &\geq F(T_{L_k}({\boldsymbol y}_k)) + (1 - \alpha_k)r_{L_k}^T ({\boldsymbol y}_k) \left( {\boldsymbol x}_k - {\boldsymbol y}_k \right) + \sum_{j = 1}^{k} \frac{\beta_{j, k+1} \gamma_j}{2} \| {\boldsymbol x}_{\Phi_{k+1}}^* - {\boldsymbol v}_j \|^2 + \left( \frac{1}{2L_k} - \frac{\alpha_k^2}{2\gamma_{k+1}} \right) \| r_{L_k} ({\boldsymbol y}_k) \|^2 \nonumber \\ 
		&+ \frac{\alpha_k \gamma_k (1 - \alpha_k)}{\gamma_{k+1}} ({\boldsymbol v}_k - {\boldsymbol y}_k)^T r_{L_k} ({\boldsymbol y}_k) + \frac{\alpha_k^2 (1 - \alpha_k) }{\gamma_{k+1}} \sum_{j = 1}^{k-1} \beta_{j, k} \gamma_j ({\boldsymbol v}_j - {\boldsymbol y}_k)^T r_{L_k} ({\boldsymbol y}_k). \label{38}
	\end{align}
\end{figure*}	
\begin{figure*}
	\setcounter{equation}{32}
	\begin{align}
		\Phi_{k+1}^*  &\geq \! F(T_{L_k}({\boldsymbol y}_k)) + (1 - \alpha_k) r_{L_k}^T ({\boldsymbol y}_k) \left( {\boldsymbol x}_k - {\boldsymbol y}_k \right) + \left( \frac{1}{2L_k} - \frac{\alpha_k^2}{2\gamma_{k+1}} \right) \| r_{L_k} ({\boldsymbol y}_k) \|^2  + \frac{\alpha_k \gamma_k (1 - \alpha_k)}{\gamma_{k+1}} ({\boldsymbol v}_k - {\boldsymbol y}_k)^T r_{L_k} ({\boldsymbol y}_k) \nonumber \\ 
		&+ \frac{\alpha_k^2 (1 - \alpha_k)}{\gamma_{k+1}} \sum_{j = 1}^{k-1} \beta_{j, k} \gamma_j ({\boldsymbol v}_j - {\boldsymbol y}_k)^T r_{L_k} ({\boldsymbol y}_k). \label{39}
	\end{align}
\end{figure*}	
Adding $\frac{\gamma_{k+1}}{2} \| {\boldsymbol x}_{\Phi_{k+1}}^* - {\boldsymbol v}_{k+1} \|^2$ to the left-hand side (LHS) of \eqref{38}, as well as moving the term $\sum_{i = 1}^{k} \frac{\beta_{i, k+1} \gamma_i}{2} \| {\boldsymbol x}_{\Phi_{k+1}}^* - {\boldsymbol v}_i \|^2$ to the LHS, we arrive to inequality \eqref{39} for $\Phi_{k+1}^*$ shown at the top of the next page.

From \eqref{39}, we have
\begin{align}
	\label{alpha_k_intuition}
	\alpha_k = \sqrt{\frac{\gamma_{k+1}}{L_k}}.
\end{align}
Substituting \eqref{gamma_expr} into \eqref{alpha_k_intuition}, the solution for $\alpha_k$ is found as 
\begin{align}
	\label{alpha_k_SFGM}
	\alpha_{k} = \frac{\sigma_k - \gamma_{k} + \sqrt{(\sigma_k - \gamma_{k})^2 + 4{L_k} \gamma_{k}}}{2{L_k}}.
\end{align}
This allows to simplify \eqref{39} as
	\begin{align}
		\label{sfgm_eq_9}
		\begin{split}
			\Phi_{k+1}^*  &\geq F(T_{L_k}({\boldsymbol y}_k)) + (1 - \alpha_k) r_{L_k}^T ({\boldsymbol y}_k) \left( {\boldsymbol x}_k - {\boldsymbol y}_k \right) \\ 
			&+ \frac{\alpha_k^2 (1 - \alpha_k)}{\gamma_{k+1}}\sum_{j = 1}^{k-1} \beta_{j, k} \gamma_j ( {\boldsymbol v}_j - {\boldsymbol y}_k)^T r_L ({\boldsymbol y}_k) \\ 
			&+ \frac{\alpha_k \gamma_k (1 - \alpha_k)}{\gamma_{k+1}} ({\boldsymbol v}_k - {\boldsymbol y}_k)^T r_L ({\boldsymbol y}_k).
		\end{split}
	\end{align}
	
\noindent Next, let us set
\begin{align}
	{\boldsymbol x}_k - {\boldsymbol y}_k + \frac{ \alpha_k \gamma_k }{\gamma_{k+1}} ({\boldsymbol v}_k - {\boldsymbol y}_k) + \frac{\alpha_k^2}{\gamma_{k+1}}\sum_{j = 1}^{k-1} \beta_{j, k} \gamma_j ({\boldsymbol v}_j - {\boldsymbol y}_k) = 0,
\end{align}
which yields
\begin{align}
	\label{y_k}
	{\boldsymbol y}_k &=  \frac{\gamma_{k+1} {\boldsymbol x}_k + \alpha_k \gamma_{k} {\boldsymbol v}_{k} + \alpha_k^2 \sum\limits_{j = 1}^{k-1} \beta_{j, k} \gamma_j {\boldsymbol v}_j}{\gamma_{k+1} + \alpha_k \gamma_{k} + \alpha_k^2 \sum\limits_{j = 1}^{k-1} \beta_{j, k} \gamma_j}. 
\end{align}
Letting ${\boldsymbol x}_{k+1} = T_{L_k} ({\boldsymbol y}_k)$ ensures that $\Phi_{k+1} \geq F ({\boldsymbol x}_{k+1})$. 

Before introducing our proposed method, let us also present a backtracking line-search strategy that will enable the convergence of the minimization process.\footnote{Several backtracking strategies have already been proposed in the literature (see \cite{Nesterov_2007, FISTA}, for example).} Since the true values of $L_{\hat{f}}$ and $\mu_{\hat{f}}$ are not known, and considering the typical applications \cite{Iulian_2}, we prioritize: \textit{(i)}~robustness to the imperfect initialization of the estimate of $L$ at iteration $k=0$; \textit{(ii)}~the need to adjust the value of the estimates of $L_{\hat{f}}$. This is achieved by selecting the parameters $\eta_u > 1$ and $\eta_d \in ]0,1[$, which are used to increase and decrease the estimate of $L_{\hat{f}}$ across different iterations. Considering this choice of parameters $\eta_u, \eta_d$, despite the initialization of $L_0$, we can always write
\begin{align}
	\label{L_bound}
	L_k \leq L_{\text{max}} \triangleq \text{max} \{\eta_d L_0, \eta_u L_{\hat{f}}\}.
\end{align}
We conclude by outlining our proposed method in Algorithm~\ref{FGM}. {In Algorithm~\ref{FGM}, $\hat{(\cdot)}$ refers to a computed value by the algorithm, and $K_{\text{max}}$ denotes the maximum number of iterations, which is linked to the tolerance via the inequalities derived in the following section that bound the difference $F({\boldsymbol x}_k) - F({\boldsymbol x}^*)$. In practice, the tolerance $\varepsilon$ is first set up. The algorithm is considered to converge with a given tolerance if $F({\boldsymbol x}_k) - F({\boldsymbol x}^*) \leq \varepsilon$. Then an upperbound on $K_{\text{max}}$ can be estimated using the bound on $F({\boldsymbol x}_k) - F({\boldsymbol x}^*)$. Such bound is tight for first-order methods.}

\begin{algorithm}
	\caption{Proposed Method}
	\label{FGM}
	\begin{algorithmic}[1]
		\STATE{\textbf{Input} ${\boldsymbol x}_0 \in \mathcal{R}^n$, $L_0 > 0$, $\mu_{\hat{f}}$, $\gamma_{0} \in [0, \mu_{\hat{f}} {]} \cup [2\mu_{\hat{f}}, 3L_0 + \mu_{\hat{f}}]$, \newline $\eta_u > 1$ and $\eta_d \in [0,1]$.}
		\STATE{\textbf{Set }$k = 0$, $i = 0$ and ${\boldsymbol v}_0 = {\boldsymbol x}_0$.}
		\WHILE{$k \leq K_{\text{max}}$}
		\STATE{$\hat{L}_i \leftarrow \eta_d L_k$}
		\WHILE{True}
		\STATE{$
			\hat{\alpha}_{i} \leftarrow \frac{\mu_{\hat{f}} + \sum\limits_{j = 1}^{k\!-\!1} \beta_{j, k} \hat{\gamma}_{j} \!-\! {\gamma}_{k} \!+\! \sqrt{\left(\mu_{\hat{f}} +\! \sum\limits_{j = 1}^{k\!-\!1} \beta_{j, k} \hat{\gamma}_{j} \!- {\gamma}_{k} \right)^2 \!+ 4\hat{L}_i {\gamma}_{k}}}{2\hat{L}_i}$}
		\STATE{$ \label{133} \hat{\gamma}_{i+1} \leftarrow (1 - \hat{\alpha}_i) {\gamma}_{k} + \hat{\alpha}_i \Big( \mu_{\hat{f}} + \sum\limits_{j = 1}^{k-1} \beta_{j, k} \hat{\gamma}_j \Big) $}
		\STATE{$\hat{{\boldsymbol y}}_i \leftarrow \frac{\hat{\gamma}_{i+1} {\boldsymbol x}_k + \hat{\alpha}_i {\gamma}_{k} {\boldsymbol v}_{k} + \hat{\alpha}_i^2 \sum\limits_{j = 1}^{k-1} \beta_{j, k} \hat{\gamma}_{j} {\boldsymbol v}_j}{\hat{\gamma}_{i+1} + \hat{\alpha}_i {\gamma}_{k} + \hat{\alpha}_i^2 \sum\limits_{j = 1}^{k-1} \beta_{j, k} \hat{\gamma}_{j}}$} 
		\STATE{ $\hat{{\boldsymbol x}}_{i+1} \leftarrow \text{prox}_{\frac{1}{\hat{L}_i} \hat{g}} \left( \hat{{\boldsymbol y}}_i - \frac{1}{\hat{L}_i} \nabla f(\hat{{\boldsymbol y}}_i) \right)$ }
		\STATE{$\hat{{\boldsymbol v}}_{i+1} \leftarrow \frac{1}{\hat{\gamma}_{i+1}} \bigg( \! (1 \!-\! \hat{\alpha}_i)\gamma_k {\boldsymbol v}_{k} \! + \! \hat{\alpha}_i \Big(\!\mu_{\hat{f}} \hat{{\boldsymbol y}}_i \!+\! \sum\limits_{j = 1}^{k-1} \beta_{j, k} \hat{\gamma}_j \hat{{\boldsymbol v}}_j $ \\ $ \qquad \qquad - \hat{L}_i \left(\hat{{\boldsymbol y}}_i \! - \! \hat{{\boldsymbol x}}_{i+1} \right) \! \Big) \! \bigg)$}
		\IF {$F (\hat{{\boldsymbol x}}_{i+1}) \leq m_{\hat{L}_i} (\hat{{\boldsymbol y}}_i, \hat{{\boldsymbol x}}_{i+1})$} 
		\STATE Break from loop
		\ELSE
		\STATE $\hat{L}_{i+1} \leftarrow \eta_u \hat{L}_i$
		\ENDIF 
		\STATE{$i \leftarrow i+1$}
		\ENDWHILE
		\STATE{$L_{k+1} \leftarrow \hat{L}_i$, ${\boldsymbol x}_{k+1} \leftarrow \hat{{\boldsymbol x}}_{i}$, $\alpha_k \leftarrow \hat{\alpha}_{i-1}$, ${\boldsymbol y}_k \leftarrow \hat{{\boldsymbol y}}_{i-1}$, $\gamma_{k+1} \leftarrow \hat{\gamma}_i$, ${\boldsymbol v}_{k+1} \leftarrow \hat{{\boldsymbol v}}_i$,  $i \leftarrow 0$, $k \leftarrow k+1$}
		\ENDWHILE
		\STATE{\textbf{Output} ${\boldsymbol x}_k$}
	\end{algorithmic}
\end{algorithm}

Comparing our proposed method to FGM, we can observe (from lines 6 and 7 in Algorithm \ref{FGM}) the differences in computing the iterates $\alpha_k$ and $\gamma_k$. In our case, their values are also dependent on the memory terms that were used in devising the estimating sequences. The update of ${\boldsymbol y}_k$ is also different, and independent of $\mu_{\hat{f}}$. A major difference is the update for ${\boldsymbol x}_k$, which is now done through a proximal gradient step. The last difference between the methods can be observed from the update of the iterates ${\boldsymbol v}_k$, which now depend on the selected subgradient. Further, comparing our proposed method to the one presented in \cite{Dosti_1} for minimizing convex functions with composite structure, we can see that the major differences arise from making use of the additional memory terms. Note that our proposed method reduces: a) to FGM when $\tau = 0$ and $\psi_k ({\boldsymbol x}) = 0, \, k = 0, 1, \ldots$, and b) to the method presented in \cite{Dosti_1} when $\psi_k ({\boldsymbol x}) = 0, \, k = 0, 1, \ldots$. In this sense, our proposed method is a generalization of all the aforementioned estimating sequence methods.

\section{Convergence Analysis}
\label{sec: Convergence Analysis}
Based on Lemma~\ref{SFGM_lemma_1}, the convergence rate of the minimization process is controlled by the rate at which the terms $\{\lambda_k\}_k$ decrease and the rate at which the terms $\{\psi_k\}_k$ increase. 

\begin{theorem}
	\label{conv_analysis_t_1}
	Let $\lambda_0 = 1$ and $\lambda_k = \prod_{j = 0}^{k-1} \left(1 - \alpha_j \right)$. Then Algorithm~\ref{FGM} generates a sequence of points $\{{\boldsymbol x}_k\}_{k}$ such that
	\begin{align}
		F({\boldsymbol x}_k) - F({\boldsymbol x}^*) &\leq \lambda_{k} \left( F({\boldsymbol x}_0) - F({\boldsymbol x}^*) + \frac{\gamma_{0}}{2} \| {\boldsymbol x}_0 - {\boldsymbol x}^* \|^2 \right) \nonumber \\
		&- (1 - \lambda_k) \psi_k ({\boldsymbol x}) .
	\end{align} 
\end{theorem}
For the proof see Appendix {\it Proof of Theorem~\ref{conv_analysis_t_1}.}

Let us now establish the rate at which the terms $\{\lambda_k\}_k$ decrease.
\begin{lemma}
	\label{conv_analysis_lemma_1}
	For all $k \geq 0$, Algorithm~\ref{FGM} guarantees that 
	\begin{enumerate}
		\item If $\gamma_0 \in [0, \mu_{\hat{f}}[$, then
			\begin{align}
				\label{FGM_conv_eq_66} 
				\lambda_{k} &\leq \frac{2 \mu_{\hat{f}}}{L_k \left(e^{\frac{k + 1}{2} \sqrt{\frac{\sigma_k}{L_k}}} \! - \! e^{-\frac{k + 1}{2} \sqrt{\frac{\sigma_k}{L_k}}} \right)^2} \leq \frac{2}{(k+1)^2}.
			\end{align}
		\item If $\gamma_0 \in [2\mu_{\hat{f}}, 3L_0 + \mu_{\hat{f}}]$, then 
		\begin{align}
			\lambda_{k} \! &\leq \! \frac{4 \mu_{\hat{f}}}{(\gamma_0 - \mu_{\hat{f}}) \! \left(e^{\frac{k + 1}{2} \sqrt{\frac{\sigma_k}{L_k}}} \! - \! e^{-\frac{k + 1}{2} \sqrt{\frac{\sigma_k}{L_k}}}\right)^2} \nonumber \\
			&\leq \frac{4L_{k}}{(\gamma_0 - \mu_{\hat{f}}) (k+1)^2}.
			\label{FGM_conv_eq_66_second_interval} 
		\end{align}
	\end{enumerate}
\end{lemma}
For the proof see Appendix {\it Proof of Lemma~\ref{conv_analysis_lemma_1}.}

Compared to \cite[Lemma~2.2.4]{Nesterov_book}, Lemma~\ref{conv_analysis_lemma_1} exhibits the following benefits: \textit{(i)} Convergence of our proposed method is established also for the cases when the exact value of $L_{\hat{f}}$ is not known. \textit{(ii)} Our proposed method converges for a broader range of $\gamma_0$. Such a result is relevant because it enables the robustness of the initialization of our proposed method in the absence of the true value of $\mu_{\hat{f}}$. 

Finally, the accelerated convergence rate for the proposed method is given by the following theorem.
\begin{theorem}
	\label{th3}
	Algorithm \ref{FGM} generates a sequence of points such that
	\begin{enumerate}
		\item If $\gamma_0 \in [0, \mu_{\hat{f}}]$, then
			\begin{align}
				\label{FGM_conv_eq_66_} 
				F(x_k) \! - \! F(x^*) &\leq \frac{\mu_{\hat{f}} (L_0 + \gamma_0) \| {\boldsymbol x}_0 - {\boldsymbol x}^* \|^2}{ L_k \left( e^{\frac{k + 1}{2} \sqrt{\frac{\sigma_k}{L_k}}} \!-\! e^{-\frac{k + 1}{2} \sqrt{\frac{\sigma_k}{L_k}}} \right)^2} 
			\end{align}
		\item If $\gamma_0 \in [2\mu_{\hat{f}}, 3 L_0 + \mu_{\hat{f}}]$, then
			\begin{align}
				\label{FGM_conv_eq_66_second_interval_} 
				F(x_k) \! - \! F(x^*) \! &\leq \! \frac{2 \mu_{\hat{f}} (L_0 + \gamma_0) \| {\boldsymbol x}_0 - {\boldsymbol x}^* \|^2 }{(\gamma_0 - \mu_{\hat{f}}) \left( e^{\! \frac{k + 1}{2} \sqrt{\frac{\sigma_k}{L_k}}} \! - \! e^{\! -\frac{k + 1}{2} \sqrt{\frac{\sigma_k}{L_k}}} \right)^2} .
			\end{align} 
	\end{enumerate}
\end{theorem}
For the proof see Appendix {\it Proof of Theorem~\ref{th3}.}

{It is worth noticing the differences between the above theorem and Theorem~3 in \cite{Dosti_1}. The structural similarity between the theorems is the consequence of the fact that the same estimating function that was introduced in \cite{Dosti_1} and used there for constructing the corresponding composite estimating sequences is also used here, but for constructing the generalized composite estimating sequences extended by the term $\psi_k ({\boldsymbol x})$ (see \eqref{def_1}). The fundamental difference of using the generalized composite estimating sequences for constructing the accelerated optimization algorithm appears in the expression for the multiplicative constant for the linear convergence rate in Theorem~\ref{th3}. It now depends on $\sigma_k$ given by \eqref{sigmaK}, instead of $\mu_{\hat{f}}$ as in Theorem~3 in \cite{Dosti_1}. In turn, $\sigma_k$ depends on both  $\mu_{\hat{f}}$ and $\sum_{j = 1}^{k-1} \beta_{j, k} \gamma_j$. Thus, $\sigma_k \geq \mu_{\hat{f}}$, and the bound for the difference between $\sigma_k$ and $\mu_{\hat{f}}$ is given by inequality \eqref{psi_bound_}, which in turn defines how much the multiplicative constant 
of the linear convergence rate is guaranteed to improve compared to that in Theorem~3 of \cite{Dosti_1}.}

{It is also worth noticing that we aim to provide here a measure of the convergence rate for the proposed algorithm in the challenging framework of the unknown Lipschitz constant with an account of the memory of the algorithm through the memory term in our generalized estimating sequences. The above results, however, are applicable also to the case when the Lipschitz constant is known/estimated and fixed. In this sense, our convergence results generalize the existing convergence results typically derived for known Lipschitz constant. The particular dynamics of the change of $\sigma_k / L_k$ is driven by the backtracking procedure for $L_k$ and the term $\sum_{j = 1}^{k-1} \beta_{j, k} \gamma_j$ that comes from the memory term in the generalized estimating sequences. Describing such a dynamic analytically appears to be hard, if feasible. Thus, in our further developments, we rely on numerical studies based on simulated and real-world data. 
Intuitively, despite the presence of the term $\sigma_k / L_k$ in the multiplicative constant of the convergence rate expressions, the result is that the rate is linear for strongly convex functions. The presence of such term affects the slope of the convergence curve. 
With backtracking for $L_k$, the slope of the convergence curve is expected to be steeper than for fixed $L$, especially if $L$ is overestimated, because backtracking helps to improve the condition number estimate at each iteration. This is in line with the other first-order algorithms extended with the backtracking procedure such as, for example, Algorithm~20 in \cite{d2021acceleration} (see Corollary~4.23 there) and Algorithm~2 in \cite{Iulian_1}. Using the bound \eqref{psi_bound_} for $\sum_{j = 1}^{k-1} \beta_{j, k} \gamma_j$, we also conclude that the slope of the convergence curve for the proposed algorithm should be steeper, when the strong convexity parameter is larger, which we next investigate in terms of numerical studies.}

\section{Numerical Studies}
\label{sec: Simulations}
We now present the numerical performance of our proposed method and compare it to the existing black-box benchmarks, specifically AMGS and FISTA. We consider both quadratic and logistic loss functions. To simulate very ill-conditioned instances of our selected problems, we also use an elastic net regularizer and select different values of the hyperparameters. Throughout all the tested instances, we demonstrate the efficiency of our proposed method when compared to the selected benchmarks. In our simulations, we make use of both synthetic and real-world datasets, the latter being chosen from the Library for Support Vector Machines \cite{LIBSVM}. Moreover, throughout our simulations, we find ${\boldsymbol x}^*$ by using CVX \cite{CVX}.

We choose the terms $\beta_{j, k} = \min \left(1, \frac{\mu}{\gamma_{k-1}} \right)$, for $j = k-1$. Depending on the selection of the terms $\gamma_0$, we will consider the following instances of our proposed method: \textit{(i)}~We set $\gamma_{0} = 0$, and refer to it as ``Proposed 1''; \textit{(ii)}~We set $\gamma_{0} = \mu_{\hat{f}}$, refer to it as ``Proposed 2''; \textit{(iii)}~We set $\gamma_0 = 3 L_0 + \mu_{\hat{f}}$, and refer to it as ``Proposed 3''. To estimate the value of the Lipschitz constant for AMGS and FISTA, we make use of the line-search strategies introduced in the corresponding papers \cite{Nesterov_2007, FISTA}. Last, in all the computational examples shown below, we select the point ${\boldsymbol x}_0$ at random and use it as a starting point for all the algorithms that are compared. 
\begin{figure*}
	\begin{minipage}[b]{.48\linewidth}
		\centering
		\centerline{\includegraphics[width=1 \columnwidth]{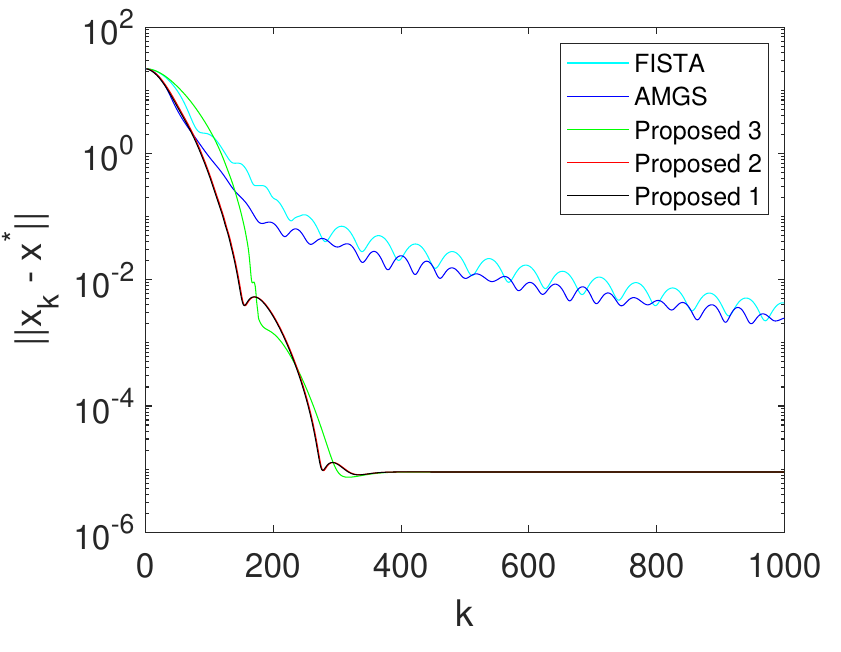}}
		\centerline{(a)}\medskip 
	\end{minipage}
	\hfill
	\begin{minipage}[b]{0.48\linewidth}
		\centering
		\centerline{\includegraphics[width=8.8cm,height=6.6cm]{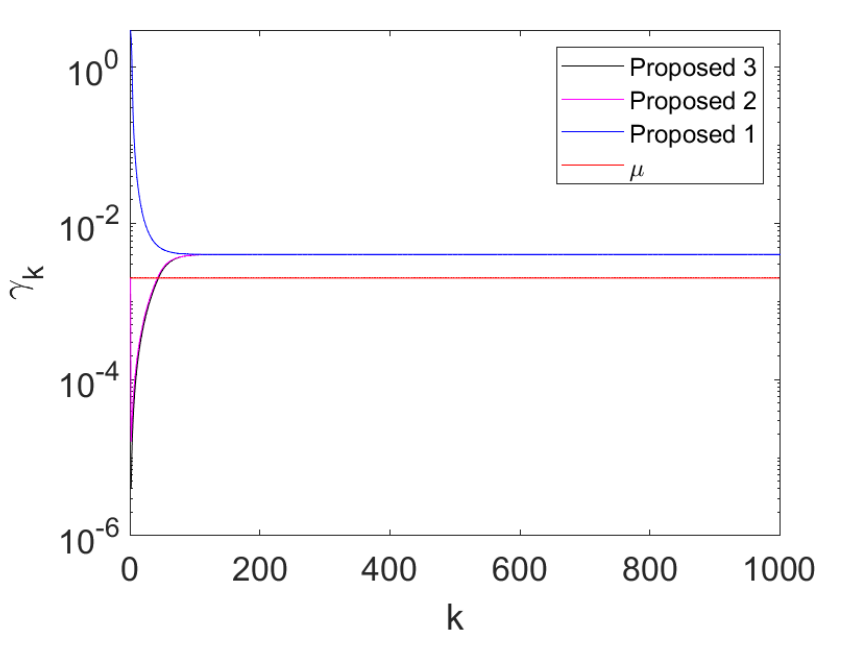}}
		\centerline{(b)}\medskip  
	\end{minipage}
	\hfill
	\begin{minipage}[b]{0.48\linewidth}
		\centering
		\centerline{\includegraphics[width=8.8cm,height=6.6cm]{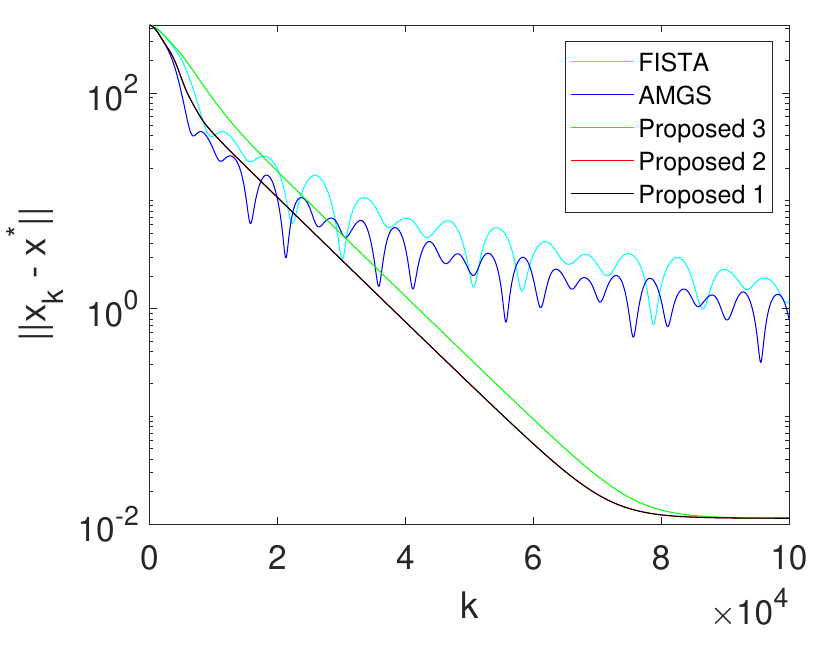}}
		\centerline{(c)}\medskip  
	\end{minipage}
	\hfill
	\begin{minipage}[b]{0.48\linewidth}
		\centering
		\centerline{\includegraphics[width=8.8cm,height=6.6cm]{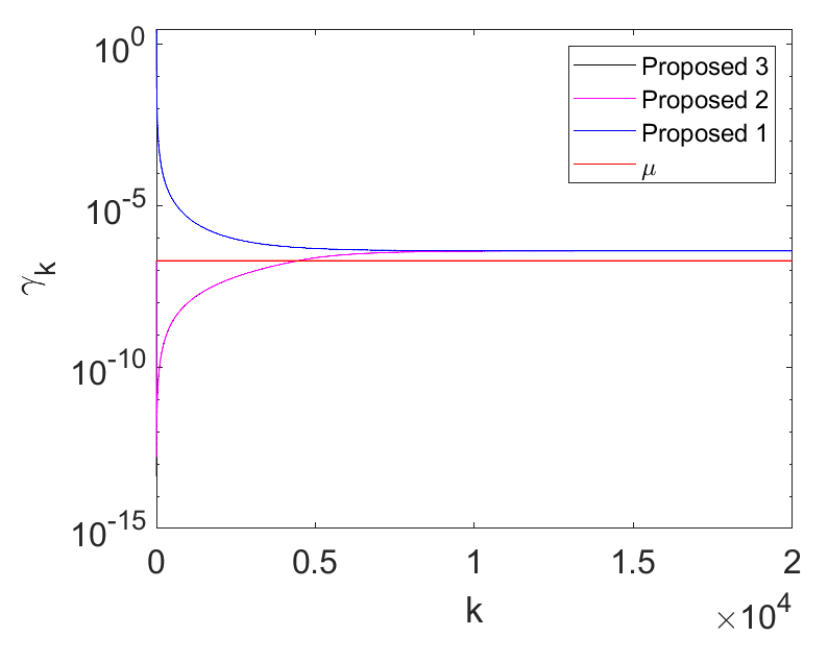}}
		\centerline{(d)}\medskip  
	\end{minipage}
	\caption{Performance evaluation of our proposed method and the selected benchmarks on synthetic data. We consider quadratic objective function and elastic net regularizer. (a)~Evaluating the distance to ${\boldsymbol x}^*$, $m = 500$, $\kappa = 10^3$ and $\tau_1 = \tau_2 = 10^{-3}$. {Note that the curves for Proposed~1 and Proposed~2 almost fully overlap.} (b)~Convergence of the terms $\{\gamma_k\}_k$, $m = 500$, $\kappa = 10^3$ and $\tau_1 = \tau_2 = 10^{-3}$. {Note that the curves for Proposed~2 and Proposed~3 almost fully overlap $\forall k$, and then also fully overlap with the curve for Proposed~1 for $k$ larger than 180.} (c)~Evaluating the distance to ${\boldsymbol x}^*$, $m = 1000$, $\kappa = 10^{7}$ and $\tau_1 = \tau_2 = 10^{-7}$. {Note that the curves for Proposed~1 and Proposed~2 fully overlap, that is, Proposed 1 and Proposed 1 have completely identical performance.} (d)~Convergence of the terms $\{\gamma_k\}_k$, $m = 1000$, $\kappa = 10^7$ and $\tau_1 = \tau_2 = 10^{-7}$. {Note that the curves for Proposed~2 and Proposed~3 fully overlap $\forall k$, and then also fully overlap with the curve for Proposed~1 for $k$ larger than about 8000.}}
	\label{num_sec_fig_1}
\end{figure*}

\subsection{Minimizing Quadratic Loss Function}
\label{linear_regression}
Let us begin with the following cost function
\begin{equation}
	\begin{aligned}
		\label{mse}
		& \underset{{\boldsymbol x} \in \mathcal{R}^n}{\text{minimize}}
		& &\frac{1}{2} \sum\limits_{i=1}^m ({\boldsymbol a}_i^T {\boldsymbol x} - {b}_i)^2 + \frac{\tau_1}{2} \| {\boldsymbol x} \|^2 + \tau_2 \| {\boldsymbol x} \|_1,
	\end{aligned}
\end{equation}
where $\| \cdot \|_1$ is the $l_1$ norm. The aim is to validate our theoretical results and demonstrate that such gains are also sustained when considering the practical deployments of the proposed method. For this purpose, we thoroughly evaluate the performance of the different benchmarks with respect to different values of the condition number of the problem. In our computational analysis, we also consider cases wherein the value of the Lipschitz constant is unknown and needs to be estimated. 

Let us start our evaluations by considering the cases where the Lipschitz constant and strong convexity parameters are known. This corresponds to the simplest case to analyze and facilitates an unbiased evaluation of the efficiency of the methods that are being compared. For this setup, we will utilize simulated data which are generated by uniformly sampling $m$ elements from the set $\{10^0, 10^{-1}, 10^{-2}, \ldots, 10^{-\xi} \}$. These elements are then used to populate the diagonal of a sparse matrix ${\boldsymbol A} = [{\boldsymbol a}_1, \cdots, {\boldsymbol a}_m]\in \mathcal{R}^{m \times m}$. The other entries of ${\boldsymbol A}$ are set to $0$. Considering the design of the matrix ${\boldsymbol A}$, we have $L = 1$ and $\mu_f = 10^{-\xi}$. Thus, the condition number of the problem becomes $\kappa = 10^{\xi}$. The entries of ${\boldsymbol y} \in \mathcal{R}^m$ are uniformly sampled from the interval $[0, 1]^n$. The other simulation parameters are set to $m \in \{500, 1000\}$, $\xi \in \{3, 7\}$ and $\tau_1 = \tau_2 \in \{10^{-3}, 10^{-7} \}$.
\begin{figure*}
	\begin{minipage}[b]{.48\linewidth}
		\centering
		\centerline{\includegraphics[width=1 \columnwidth]{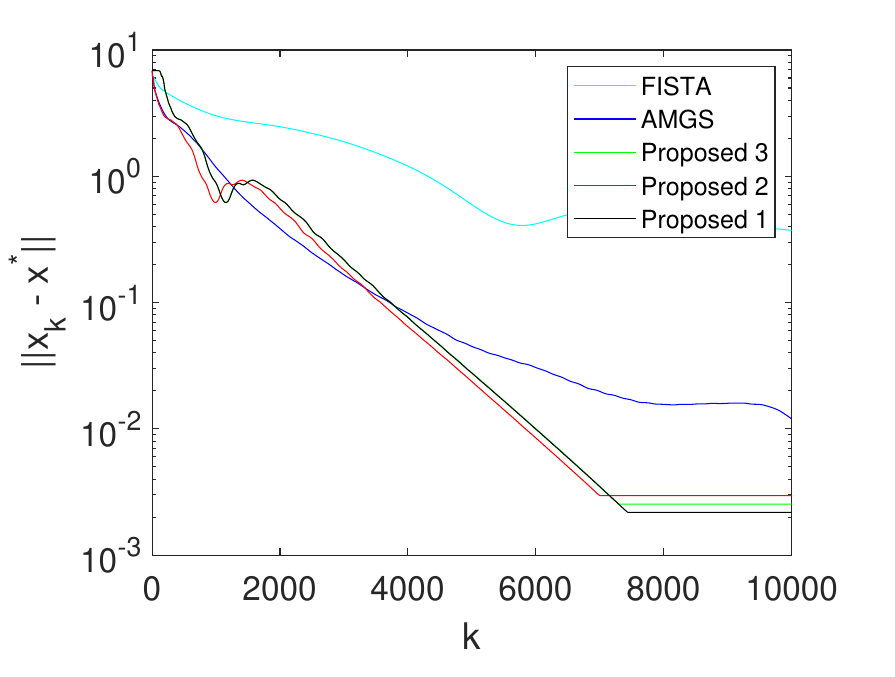}}
		\centerline{(a)}\medskip 
	\end{minipage}
	\hfill
	\begin{minipage}[b]{0.48\linewidth}
		\centering
		\centerline{\includegraphics[width=8.8cm,height=6.6cm]{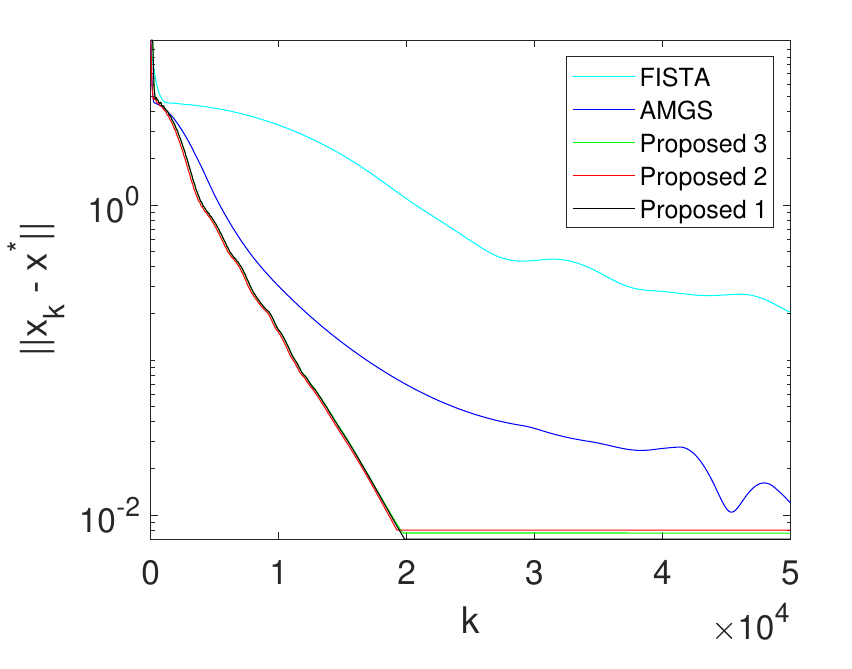}}
		\centerline{(b)}\medskip  
	\end{minipage}
	\hfill
	\begin{minipage}[b]{0.48\linewidth}
		\centering
		\centerline{\includegraphics[width=8.8cm,height=6.6cm]{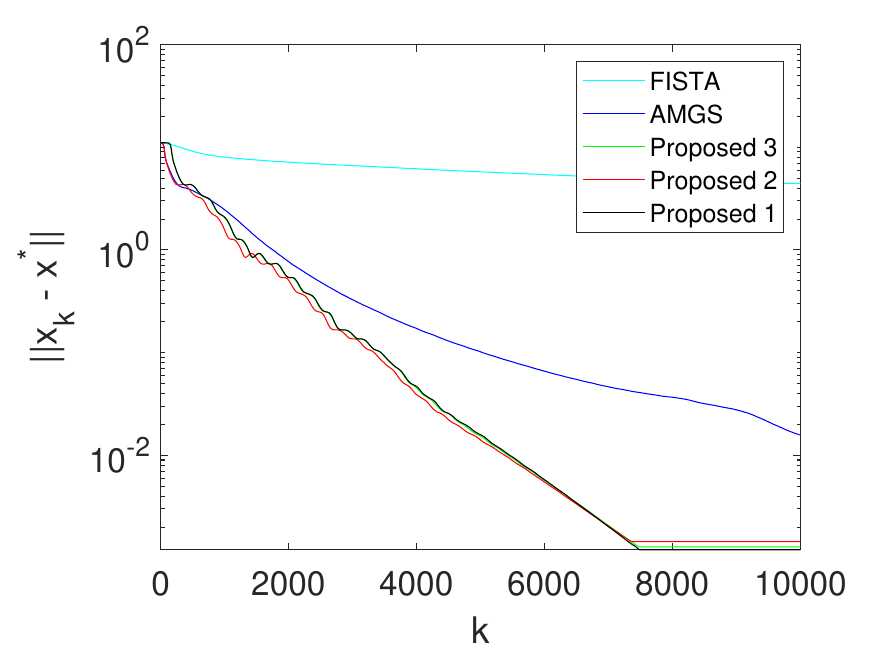}}
		\centerline{(c)}\medskip  
	\end{minipage}
	\hfill
	\begin{minipage}[b]{0.48\linewidth}
		\centering
		\centerline{\includegraphics[width=8.8cm,height=6.6cm]{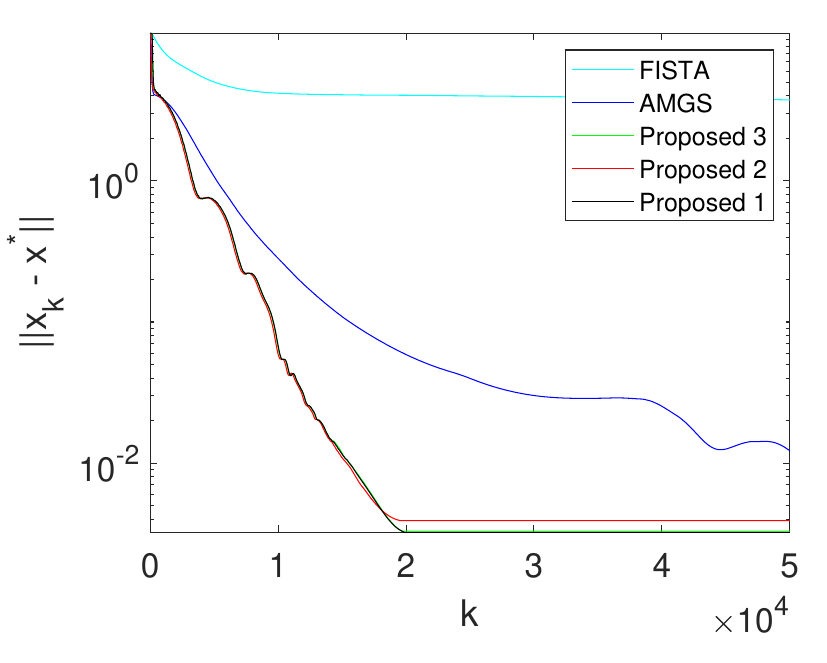}}
		\centerline{(d)}\medskip  
	\end{minipage}
	\caption{Performance evaluation of our proposed method and the selected benchmarks on the ``a1a'' dataset. We consider quadratic objective function and elastic net regularizer, and assume that the true value of $L_{\hat{f}}$ is not known. (a)~Evaluating the distance to ${\boldsymbol x}^*$ for ``a1a'' dataset, $L_0 = 0.1 L_{\text{``a1a''}}$ and $\tau_1 = \tau_2 = 10^{-4}$. (b)~Evaluating the distance to ${\boldsymbol x}^*$ for ``a1a'' dataset, $L_0 = 0.1 L_{\text{``a1a''}}$ and $\tau_1 = \tau_2 = 10^{-5}$. {Note that the curves for Proposed~2 and Proposed~3 almost fully overlap.} (c)~Evaluating the distance to ${\boldsymbol x}^*$ for ``a1a'' dataset, $L_0 = 10 L_{\text{``a1a''}}$ and $\tau_1 = \tau_2 = 10^{-4}$. (d)~Evaluating the distance to ${\boldsymbol x}^*$ for ``a1a'' dataset, $L_0 = 10 L_{\text{``a1a''}}$ and $\tau_1 = \tau_2 = 10^{-5}$.}
	\label{num_sec_fig_2}
\end{figure*}

\begin{figure*}
	\begin{minipage}[b]{.48\linewidth}
		\centering
		\centerline{\includegraphics[width=1 \columnwidth]{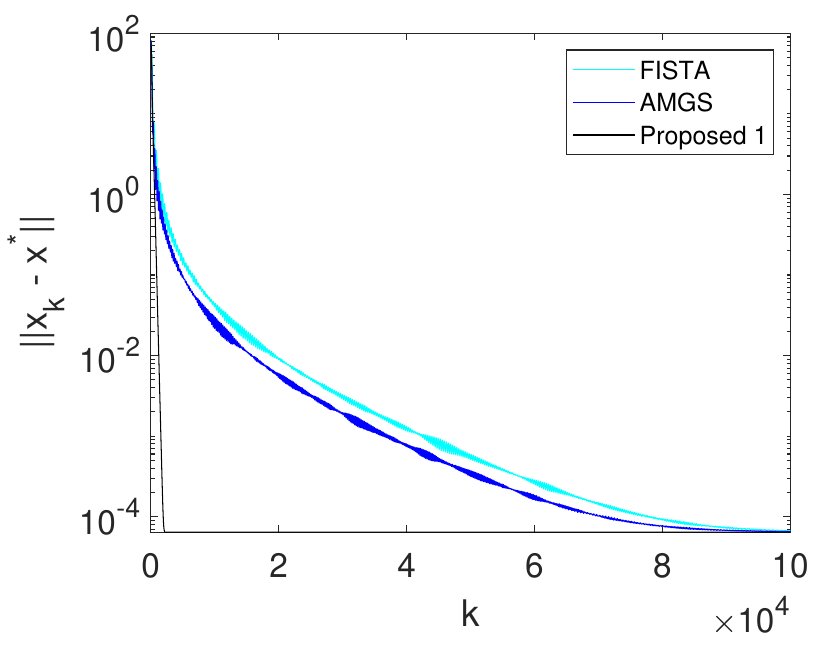}}
		\centerline{(a)}\medskip 
	\end{minipage}
	\hfill
	\begin{minipage}[b]{0.48\linewidth}
		\centering
		\centerline{\includegraphics[width=8.8cm,height=6.6cm]{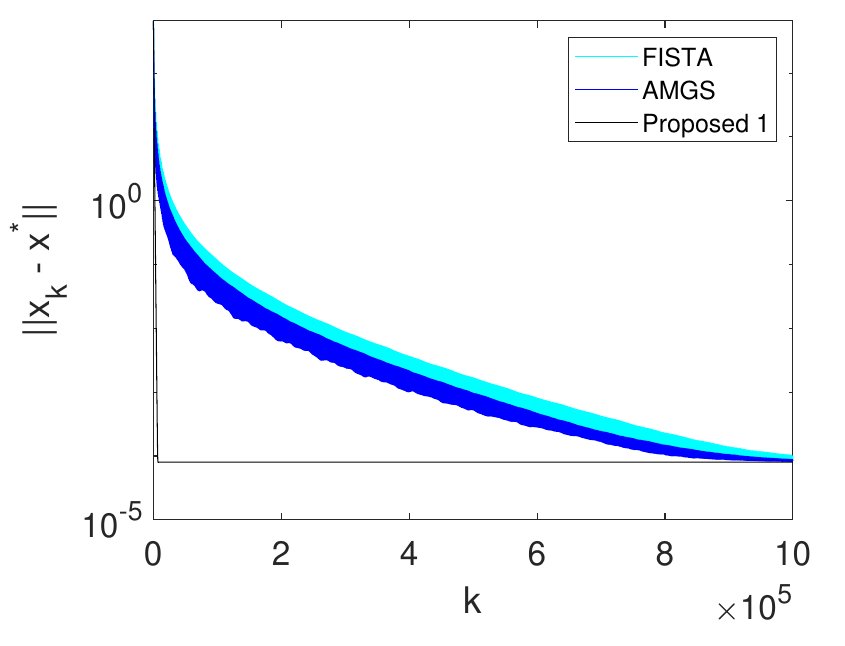}}
		\centerline{(b)}\medskip  
	\end{minipage}
	\hfill
	\begin{minipage}[b]{0.48\linewidth}
		\centering
		\centerline{\includegraphics[width=8.8cm,height=6.6cm]{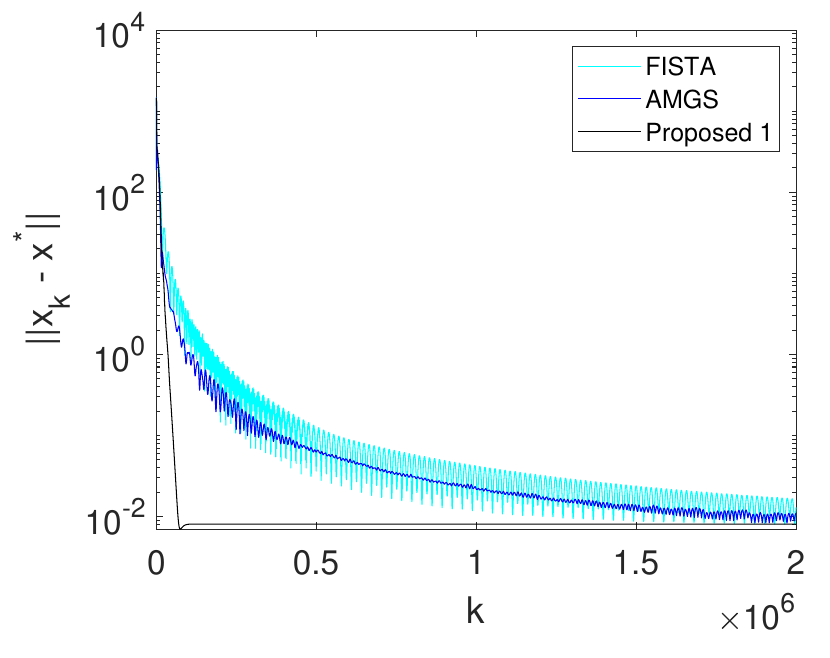}}
		\centerline{(c)}\medskip  
	\end{minipage}
	\hfill
	\begin{minipage}[b]{0.48\linewidth}
		\centering
		\centerline{\includegraphics[width=8.8cm,height=6.6cm]{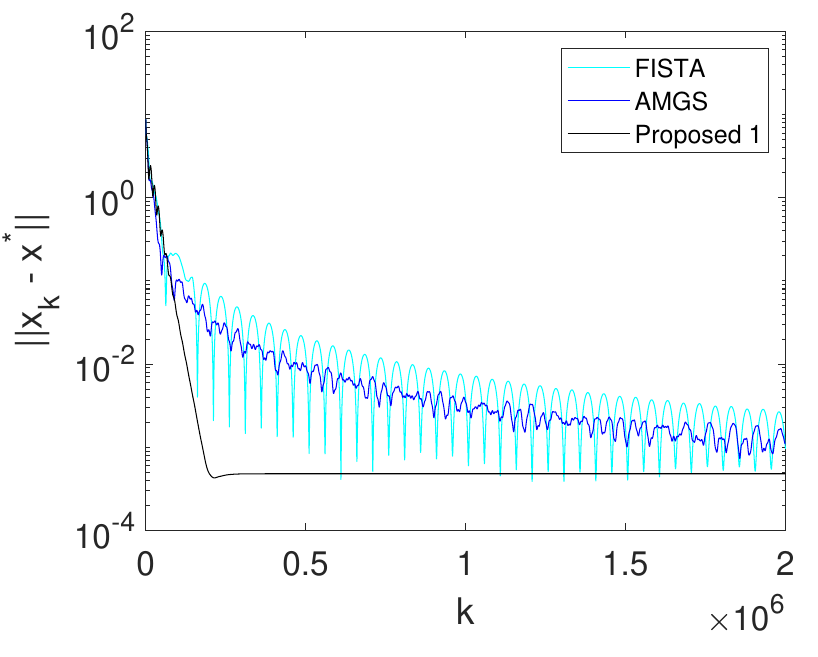}}
		\centerline{(d)}\medskip  
	\end{minipage}
	\caption{Performance evaluation of our proposed method and the selected benchmarks on real data. We consider the logistic objective function and elastic net regularizer. (a)~Evaluating the distance to ${\boldsymbol x}^*$ for ``rcv1.binary'' dataset, $\tau_1 = \tau_2 = 10^{-4}$. (b)~Evaluating the distance to ${\boldsymbol x}^*$ for ``rcv1.binary'' dataset, $\tau_1 = \tau_2 = 10^{-5}$. (c)~Evaluating the distance to ${\boldsymbol x}^*$ for ``triazine'' dataset, $\tau_1 = \tau_2 = 10^{-6}$. (d)~Evaluating the distance to ${\boldsymbol x}^*$ for ``triazine'' dataset, $\tau_1 = \tau_2 = 10^{-7}$.}
	\label{num_sec_fig_3}
\end{figure*}

When compared to the selected benchmarks, we can observe in Fig.~\ref{num_sec_fig_1} that our proposed method is more efficient both in terms of the obtained distance to the optimal solution ${\boldsymbol x}^*$, as well as in the number of iterations needed to converge to such a solution. Another advantage of our proposed method is that it exhibits better monotonic properties. Moreover, all the methods that are being evaluated are sensitive to the condition number of the problem. The higher the value of the condition number is, the more iterations the methods require to converge in the vicinity of ${\boldsymbol x}^*$. Comparing between the selected instances of our proposed method, we can observe that they exhibit a commensurate degree of similarity, which is also clear based on our theoretical analysis. Nevertheless, we can see that the best performing instance is the one obtained when choosing $\gamma_0 = 0$.

Let us next consider the case where the true value of the Lipschitz constant is not known. For this purpose, we shall consider initial estimates of the Lipschitz constant that are $10$ times higher and lower than the true value, i.e., $L_0 \in \{0.1L_{\hat{f}}, 10L_{\hat{f}}\}$. {Following the recommendations presented in \cite{becker2011templates}, for our line-search procedure we choose $\eta_u = 2$ and $\eta_d = 0.9$.} We also assume the true value of the strong convexity parameter $\mu_{\hat{f}}$ is not known. Instead, we use the lower bound on the true value, which can be controlled by the selection of the regularizer term in \eqref{mse}. In the following examples, we will use data from the fluorescent protein database ``a1a'' \cite{LIBSVM}, for which ${\boldsymbol A} \in \mathcal{R}^{1605 \times 123}$. For the considered dataset, the true value of the Lipschitz constant is $L_{\text{``a1a''}} = 10061$. The values of the regularizers are selected to be $\tau_1 = \tau_2 \in \{10^{-4}, 10^{-5}\}$, which ensures that the condition number of the problem $\kappa = \frac{L_{\hat{f}}}{\mu_{\hat{f}}}$ has a high value. 

We can observe in Fig.~\ref{num_sec_fig_2} that our proposed method is more efficient than the selected benchmark. Similar to the results presented in Fig.~\ref{num_sec_fig_1}, the iterates produced from our proposed method exhibit better monotonic properties and have the smallest distance to the optimal solution. Moreover, across all simulations, we can observe that our proposed method converges to ${\boldsymbol x}^*$ in a smaller number of iterations. Considering the result for different values of regularizers and Lipschitz constant estimates, we can observe the robustness of our proposed method and AMGS to the imperfect selection of $L_0$. A difference between these two methods, however, is that AMGS exhibits a higher per-iteration complexity. Such results cannot be observed for FISTA, whose performance is very sensitive to the initialization of the Lipschitz constant estimate. This comes because the line-search strategy introduced for FISTA does not allow for decreasing the estimate of the Lipschitz constant across iterates. Comparing the different versions of our proposed method, we can observe that in most cases, they are equally efficient. Nevertheless, the variant obtained when initializing $\gamma_0 = 0$ is preferred because it enables the robustness of the initialization of our proposed method with respect to the imperfect knowledge of $\mu_{\hat{f}}$.


\subsection{Minimizing Logistic Loss Function}
\label{Decreasing the norm of the gradient}
We also test the performance of our algorithm and selected benchmarks in minimizing the following function
\begin{equation}
	\begin{aligned}
		\label {num_eq_2}
		& \underset{{\boldsymbol x} \in \mathcal{R}^n}{\text{minimize}}
		& &\frac{1}{m} \sum\limits_{i=1}^m \text{log} \left( 1 \!-\! \text{e}^{-b_i {\boldsymbol x}^T {\boldsymbol a}_i} \right) \!+\! \frac{\tau_1}{2} \| {\boldsymbol x} \|^2 \!+\! \tau_2 \| {\boldsymbol x} \|_1.  
	\end{aligned}
\end{equation} 

\noindent We consider datasets namely ``rcv1.binary'', for which ${\boldsymbol A}_{\text{``rcv1.binary''}} \in \mathcal{R}^{1000 \times 2000}$, and a subset of ``triazine'', for which ${\boldsymbol A}_{\text{``triazine''}} \in  \mathcal{R}^{186 \times 61}$ \cite{LIBSVM}. Moreover, we observed in the previous subsection that the convergence of FISTA is significantly affected by the selection of $L_0$, which happens because the line-search strategy proposed for FISTA does not allow for decreasing the estimate of the Lipschitz constant. Since in this paper the goal is to devise more efficient black-box algorithms, we assume that the true value of $L_{\hat{f}}$ is known. For the selected datasets, we have $L_{\text{``rcv1.binary''}} = 1.13$ and $L_{\text{``triazine''}} = 25.15$. Regarding the strong convexity parameter, we follow a similar approach as in the earlier examples and select its value to be the same as the $l_2$ regularizer term in \eqref{num_eq_2}, which are selected to be $\tau_1 = \tau_2 \in \{10^{-4}, 10^{-5}, 10^{-6}, 10^{-7}\}$. Last, since there is little performance difference between the different variants of our proposed method, in the sequel, we simulate only the first variant, namely Proposed 1. Our findings are depicted in Fig. \ref{num_sec_fig_3}, and from it, we can clearly see that our proposed method significantly outperforms the selected benchmarks in minimizing the regularized logistic loss function.


\section{Conclusion and Discussion}
\label{Con}
{A new class of estimating sequences that is named as generalized composite estimating sequences has been introduced for minimizing convex functions with composite structure with a non-smooth term. Using this newly introduced class of estimating sequences, a new accelerated black-box first-order algorithm has been proposed. The proposed algorithm is endowed with an efficient backtracking line-search strategy and exhibits an accelerated convergence rate even when the true value of the Lipschitz constant of the objective function is not known. The convergence results presented in the paper suggest that the proposed algorithm exhibits such an accelerated convergence when $\gamma_0 \in [0, 3L + \mu_{\hat{f}}]$, i.e., the initialization of our proposed method is robust to the imperfect knowledge of the strong convexity parameter as well. From a computational viewpoint, our proposed method has been shown to outperform the existing benchmarks when tested in solving practical problems for both simulated and real-world datasets.} 

{The results presented in this paper can be extended in multiple directions. First, it would be of interest to explore other structures for $\psi_k ({\boldsymbol x})$, which can be used for devising estimating sequences applicable to different optimization methods, e.g., higher-order methods, stochastic methods, non-convex methods, etc. 
Extending the framework to the inexact oracle framework, particularly in the stochastic approximation context, is also of significant interest. Additionally, studying the impact of restarting on the practical performance of the proposed method would be valuable, although such a study is more heuristic and falls outside the scope of this paper, which focuses on developing rigorous results.
Extensions of the framework devised herein in the context of the inexact oracle framework is also of a high interest. It is especially so in the stochastic approximation framework. The study of the impact of restarting to the practical performance of our proposed method is also of interest, but such study is heuristic, and thus outside of the scope of this paper devoted to developing exact results.}

\appendix[]

\section*{Proof of Lemma~\ref{SFGM_lemma_1}}
\begin{proof}
	By the condition of Lemma~\ref{SFGM_lemma_1} and using \eqref{def_1}, we have
	\begin{align} \nonumber
		F({\boldsymbol x}_k) &\leq \Phi_{k}^* = \underset{{\boldsymbol x} \in {\mathcal{Q}}}{\min} \Phi_{k} ({\boldsymbol x}) \\
		&\stackrel{\eqref{def_1}}{\leq}  \!
		\underset{{\boldsymbol x} \in {\mathcal{Q}}}{\min} \lambda_{k} \Phi_{0}({\boldsymbol x}) + (1 - \lambda_{k}) \left( F({\boldsymbol x}) - \psi_k ({\boldsymbol x}) \right) \nonumber \\ 
		&\leq  \lambda_{k} \Phi_{0}({\boldsymbol x}^*) + (1 - \lambda_{k}) \left(F({\boldsymbol x}^*) - \psi_k ({\boldsymbol x}^*)\right).
	\end{align}
	Regrouping the terms concludes the proof. 
\end{proof}

\section*{Proof of Lemma~\ref{SFGM_lemma_2}}
\begin{proof}
	We prove this by induction. At step $k = 0$, considering \eqref{def_1} together with the facts that $\lambda_0 = 1$ and $\psi_0 ({\boldsymbol x}) = 0$, we can write: $\Phi_0 ({\boldsymbol x}) \leq \lambda_0 \Phi_0 ({\boldsymbol x}) + \left( 1 - \lambda_0 \right) F ({\boldsymbol x}) \equiv \Phi_0 ({\boldsymbol x})$. At iteration $k$, assume \eqref{def_1} holds true, which results in
	\begin{align}
		\label{useful}
		\Phi_k({\boldsymbol x}) - \left(1 - \lambda_k \right) F({\boldsymbol x}) \leq \lambda_k \Phi_0 ({\boldsymbol x}) - \left(1 - \lambda_k \right) \psi_k ({\boldsymbol x}).
	\end{align}
	
	Utilizing \eqref{13} in \eqref{phi_k+1_SFGM}, yields 
	\begin{align}
		\label{crappy}
		\Phi_{k+1} ({\boldsymbol x}) &\leq (1 - \alpha_k) \left(\Phi_{k} ({\boldsymbol x}) + \psi_k ({\boldsymbol x}) \right) + \alpha_k \left( F({\boldsymbol x}) + \psi_k ({\boldsymbol x}) \right) \nonumber \\
		&- \psi_{k+1} ({\boldsymbol x}) - \Psi.
	\end{align}
	Considering that $\Psi$ is an upper bound on $\psi_k (x)$, and adding it to the right-hand side (RHS) of \eqref{crappy}, results in 
	\begin{align}
		\label{_123_}
		\Phi_{k+1} ({\boldsymbol x}) \! &\leq \! (1 \! - \! \alpha_k) \Phi_{k} ({\boldsymbol x}) \! + \! \alpha_k F({\boldsymbol x}) \! + \! (1 - \alpha_{k}) (1 - \lambda_{k}) F({\boldsymbol x}) \nonumber \\ &- (1 - \alpha_{k}) (1 - \lambda_{k}) F({\boldsymbol x}) - \psi_{k+1} ({\boldsymbol x}).
	\end{align}
	Relaxing the RHS of \eqref{_123_}, yields
	\begin{align}
		\label{3.8}
		\Phi_{k \! + \! 1} \! ({\boldsymbol x}) &\leq \! (1 \! - \! \alpha_k) \! \left( \Phi_{k} ({\boldsymbol x}) \! - \! (1 \! - \! \lambda_{k}) F({\boldsymbol x}) \right) \nonumber  \\
		&+ \! \left(\alpha_k \! + \! (1 \! - \! \lambda_{k}) (1 \! - \! \alpha_k) \right) \! F({\boldsymbol x}) \! - \! \psi_{k \! + \! 1}({\boldsymbol x}). 
	\end{align}
	Substituting \eqref{useful} in \eqref{3.8}, results in
	\begin{align}
		\label{mmmm}
		\Phi_{k\! + \! 1} ({\boldsymbol x}) \! &\leq \! (1 \! - \! \alpha_k) \lambda_{k} \left( \Phi_{0} ({\boldsymbol x}) \! - \! (1 \! - \! \lambda_k) \psi_k ({\boldsymbol x}) \right) \nonumber \\ 
		&+ \! (1 \! - \! \lambda_{k} \! + \! \alpha_{k} \lambda_{k}) F({\boldsymbol x}) \! - \! \psi_{k \! + \! 1}({\boldsymbol x}). 
	\end{align}     
	Last, relaxing the RHS of \eqref{mmmm} and using \eqref{lambda_recursive} yields
	\begin{align}
		\Phi_{k+1} ({\boldsymbol x}) &\leq \lambda_{k+1} \Phi_{0} ({\boldsymbol x}) + (1 - \lambda_{k+1}) \left( F({\boldsymbol x})- \psi_{k+1}({\boldsymbol x}) \right). 
	\end{align}     
\end{proof}

\section*{Proof of Lemma~\ref{SFGM_lemma_3}}
\begin{proof}
	Recall that for $k = 0$, we have $\psi_0 (x) = 0$. Thus, $\nabla^2 \Phi_0 (x) = \gamma_0 {\boldsymbol I}$, where ${\boldsymbol I}$ is the identity matrix. Assume that for step $k$ we have: $\nabla^2 \Phi_k(x) = \gamma_k {\boldsymbol I} - \sum_{i = 1}^{k-1} \beta_{i, k} \gamma_i {\boldsymbol I}$. For step $k+1$, consider the following
	\begin{align}
		\label{t^3}
		\nabla^2 \Phi_{k+1} ({\boldsymbol x}) &\stackrel{\eqref{phi_k+1_SFGM}}{=} \Big( (1 - \alpha_k) \gamma_{k} +\alpha_k \sigma_k - \sum_{j = 1}^{k} \beta_{j, k} \gamma_j \Big) {\boldsymbol I}. 
	\end{align}
	Massaging \eqref{t^3} we obtain
	\begin{align}
		\label{t^2}
		\gamma_{k+1} {\boldsymbol I} = \big( (1-\alpha_k)\gamma_k + \alpha_k \sigma_k \big) {\boldsymbol I}.
	\end{align}
	Substituting \eqref{gamma_expr} into \eqref{t^2} is sufficient to establish that the quadratic cannonical structure for $\{\Phi_k\}_k$ is preserved. 
	
	Let us next focus on finding the recurrent relations for the terms $\{{\boldsymbol v}_{k}\}_{k}$. First, replacing \eqref{phi} in \eqref{phi_k+1_SFGM} and making some algebraic manipulations, results in 
	\begin{align}
		\label{finding_v}
		&\phi_{k + 1}^* \! + \! \frac{\gamma_{k + 1}}{2} \| {\boldsymbol x} \! - \! {\boldsymbol v}_{k + 1} \|^2 \! = \! (1 \! - \! \alpha_k) \left( \! \phi_{k}^* \! + \! \frac{\gamma_{k}}{2} \| {\boldsymbol x} \! - \! {\boldsymbol v}_{k} \|^2 \! \right) \nonumber \\
		&- \Psi + \alpha_{k} \Big( F\left(T_{L_k} ({\boldsymbol y}_k) \right) + \psi_{k} ({\boldsymbol x}) + \frac{1}{2L_k} \| r_{L_k} ({\boldsymbol y}_k) \|^2 \nonumber \\
		&+ r_{L_k}^T ({\boldsymbol y}_k) ({\boldsymbol x} - {\boldsymbol y}_k)  + \frac{\mu_{\hat{f}}}{2} \| {\boldsymbol x} - {\boldsymbol y}_k \|^2 \Big).
	\end{align}
	
	Observe that both sides of \eqref{finding_v} are convex in ${\boldsymbol x}$. From the first-order optimality condition we have
	\begin{align}
		&\gamma_{k + 1} ({\boldsymbol x} - {\boldsymbol v}_{k + 1}) = \gamma_k (1 - \alpha_k) ({\boldsymbol x} - {\boldsymbol v}_{k}) \nonumber \\
		&\quad + \alpha_k  \left( \mu_{\hat{f}} ({\boldsymbol x} \!-\! {\boldsymbol y}_k) \!+\! r_{L_k} ({\boldsymbol y}_k) \!+\! \sum_{j = 1}^{k-1} \beta_{j, k} \gamma_j ({\boldsymbol x} \!-\! {\boldsymbol v}_j) \right). \label{SFGM_opt_cond}
	\end{align} 
	Substituting \eqref{gamma_expr} in \eqref{SFGM_opt_cond}, and reducing the dependency on ${\boldsymbol x}$ results in
	\begin{align}
		\label{v_update}
			- \gamma_{k + 1} {\boldsymbol v}_{k + 1} &= \alpha_k \left( r_{L_k} ({\boldsymbol y}_k) - \mu_{\hat{f}}  {\boldsymbol y}_k - \sum_{j = 1}^{k-1} \beta_{j, k} \gamma_j {\boldsymbol v}_j \right) \nonumber \\
			&- (1 - \alpha_k) \gamma_k {\boldsymbol v}_{k}.
	\end{align}
	Substituting \eqref{r_l} into \eqref{v_update} yields the desired \eqref{v_value}. 
	
	Let us now focus on finding the terms $\{\phi_k^*\}_k$. A straightforward approach is to assume that there exists a sequence of estimating functions $\{\Theta_k ({\boldsymbol y}_k)\}_k$ for the sequence $\{{\boldsymbol y}_k\}_k$ that has the following structure
	\begin{align}
		\label{theta_k}
		\Theta_k ({\boldsymbol y}_k) = \theta_k^* + \frac{\gamma_k}{2} \| {\boldsymbol y}_k - {\boldsymbol v}_k \|^2 - \sum_{j = 1}^{k-1} \frac{\beta_{j, k} \gamma_j}{2} \| {\boldsymbol y}_k - {\boldsymbol v}_j \|^2 
	\end{align}
	Next, consider \eqref{phi_k+1_SFGM} with ${\boldsymbol x} = {\boldsymbol y}_k$
	\begin{align}
		\label{theta_kk}
		\Theta_{k+1} ({\boldsymbol y}_k) \!  &= \! (1 \! - \! \alpha_k) \left(\Theta_{k} ({\boldsymbol y}_k) + \psi_{k} ({\boldsymbol y}_k) \right) \! - \psi_{k+1} ({\boldsymbol y}_k) - \Psi \nonumber \\  &+  \! \alpha_{k} \! \left( \! F\left(T_{L_k} \! ({\boldsymbol y}_k) \right) \! + \psi_{k} ({\boldsymbol y}_k) + \! \frac{1}{2L_k} \|r_{L_k} \! ({\boldsymbol y}_k) \|^2 \! \right).
	\end{align}
	Substituting \eqref{psii} and \eqref{theta_k} into \eqref{theta_kk}, and relaxing the RHS, results in
	\begin{align}
		\label{theta_kkk}
		&\theta_{k+1}^* + \frac{\gamma_{k+1}}{2} \| {\boldsymbol y}_k - {\boldsymbol v}_{k+1} \|^2 \! \leq \! (1 \! - \! \alpha_k) \left(\theta_{k}^* + \frac{\gamma_{k}}{2} \| {\boldsymbol y}_k - {\boldsymbol v}_{k} \|^2 \right) \nonumber \\
		&+ \! \alpha_{k} \! \Big( \! F\left(T_{L_k} \! ({\boldsymbol y}_k) \right) \! + \! \frac{1}{2L_k} \| r_{L_k} \! ({\boldsymbol y}_k) \|^2 + \sum_{j = 1}^{k-1} \frac{\beta_{j, k} \gamma_j}{2} \| {\boldsymbol y}_k \!-\! {\boldsymbol v}_j \|^2 \! \Big).
	\end{align}
	
	Using \eqref{v_value}, we can write
	\begin{align}
		\label{bound}
		{\boldsymbol v}_{k+1} \! - \! {\boldsymbol y}_k \! &= \! \frac{1}{\gamma_{k+1}} \Big( \! (1 \! - \! \alpha_k)\gamma_k {\boldsymbol v}_k \! + \! \alpha_k \Big( \! \mu_{\hat{f}} {\boldsymbol y}_k \! - \! r_{L_k} ({\boldsymbol y}_k) \nonumber \\
		&+ \! \sum_{j = 1}^{k-1} \beta_{j, k} \gamma_j {\boldsymbol v}_j \Big) - \gamma_{k+1} {\boldsymbol y}_k \Big).
	\end{align}
	Substituting \eqref{gamma_expr} into \eqref{bound}, after some straightforward algebraic manipulations, we can rewrite \eqref{bound} as
	\begin{align}
		\label{boundd}
		{\boldsymbol v}_{k+1} &- {\boldsymbol y}_k =\frac{1}{\gamma_{k+1}} \Bigg((1-\alpha_k)\gamma_k ({\boldsymbol v}_k - {\boldsymbol y}_k) \nonumber \\ 
		&+ \alpha_k \left(\sum_{j = 1}^{k-1} \beta_{j, k} \gamma_j ({\boldsymbol v}_j - {\boldsymbol y}_k) - r_{L_k} ({\boldsymbol y}_k)\right) \Bigg).
	\end{align}
	Taking $||\cdot||^2$ of \eqref{boundd}, multiplying with $\frac{\gamma_{k+1}}{2}$, and extending the RHS, we reach \eqref{fgm_eq_1} shown at the bottom of the next page. 	
	\begin{figure*}[!b]
		\setcounter{equation}{63}
	\begin{align}
		\label{fgm_eq_1}
		\frac{\gamma_{k\!+\!1}}{2} \| {\boldsymbol v}_{k\!+\!1} \! &- \! {\boldsymbol y}_k \|^2 \! = \! \frac{(1 \! - \! \alpha_k)^2\gamma_k^2}{2 \gamma_{k\!+\!1}} \| {\boldsymbol v}_{k} \! - \! {\boldsymbol y}_k \|^2 \! + \! \frac{\alpha_{k}^2}{2 \gamma_{k\!+\!1}} \! \left( \! \|r_{L_k} ({\boldsymbol y}_k) \|^2  + \! \left\| \sum_{j = 1}^{k-1} \beta_{j, k} \! \gamma_j ({\boldsymbol v}_j \! - \! {\boldsymbol y}_k) \right\|^2 \right) \nonumber \\ 
		&-\! \frac{\alpha_k^2}{\gamma_{k\!+\!1}} \sum_{j = 1}^{k - 1} \beta_{j, k} \gamma_j ({\boldsymbol v}_j \!-\! {\boldsymbol y}_k)^T r_{L_k} ({\boldsymbol y}_k) - \! \frac{\alpha_k (1 \!- \! \alpha_k)\gamma_{k}}{\gamma_{k+1}} \left( \! ({\boldsymbol v}_{k} \! - \! {\boldsymbol y}_k)^T \! r_{L_k} \left( {\boldsymbol y}_k \right) \! - \! \sum_{j = 1}^{k-1} \! \beta_{j, k} \gamma_j ({\boldsymbol v}_j \! - \! {\boldsymbol y}_k)^T ({\boldsymbol v}_k \! - \! {\boldsymbol y}_k) \! \right) .
	\end{align}
	\end{figure*}	
	Substituting \eqref{fgm_eq_1} into \eqref{theta_kkk}, yields \eqref{theta_kkkk} shown at the bottom of the next page. 
	\begin{figure*}[!b]
		\setcounter{equation}{64}
		\begin{align}
			\label{theta_kkkk}
			\theta_{k+1}^* \!  &\leq \! (1 \! - \! \alpha_k) \theta_{k}^* + \frac{(1-\alpha_k) \gamma_k}{2} \left( \! 1 \!-\! \frac{(1-\alpha_k)\gamma_k}{\gamma_{k+1}} \! \right) \| {\boldsymbol y}_k - {\boldsymbol v}_{k} \|^2 \!+\! \alpha_{k} \! \left( \! F\left(T_{L_k} \! ({\boldsymbol y}_k) \right) \!+\! \frac{1}{2L_k} \| r_{L_k} \! ({\boldsymbol y}_k) \|^2 \!+\! \sum_{j = 1}^{k-1} \frac{\beta_{j, k} \gamma_j}{2} \| {\boldsymbol v}_j - {\boldsymbol y}_k \|^2 \right) \nonumber \\
			&- \frac{\alpha_k^2}{2\gamma_{k+1}} \left( \left\| \sum_{j = 1}^{k-1} \frac{\beta_{j, k} \gamma_j}{2} ({\boldsymbol y}_k - {\boldsymbol v}_j) \right\|^2 \!+\! \| r_{L_k} ({\boldsymbol y}_k) \|^2 \right) \!+\! \frac{\alpha_k^2}{\gamma_{k+1}} \sum_{j = 1}^{k-1} \beta_{j, k} \gamma_j ({\boldsymbol v}_j - {\boldsymbol y}_k)^T r_{L_k} ({\boldsymbol y}_k) \nonumber \\ 
			&+ \frac{\alpha_k (1-\alpha_k)\gamma_{k}}{\gamma_{k+1}} \left( ({\boldsymbol v}_{k} - {\boldsymbol y}_k)^T r_{L_k} \left( {\boldsymbol y}_k \right) - \sum_{j = 1}^{k-1} \beta_{j, k} \gamma_j ({\boldsymbol v}_j - {\boldsymbol y}_k)^T ({\boldsymbol v}_k - {\boldsymbol y}_k) \! \right).
		\end{align}
	\end{figure*}	
	In \eqref{theta_kkkk}, using the Cauchy-Schwartz inequality and relaxing the upper bound, yields in turns \eqref{theta_kkkkk} shown at the bottom of the next page. 
		\begin{figure*}[!b]
		\setcounter{equation}{65}
		\begin{align}
			\label{theta_kkkkk}
			\theta_{k+1}^* \!  &\leq \! (1 \! - \! \alpha_k) \theta_{k}^* + \frac{\alpha_k \gamma_k (1 \!-\! \alpha_k) \sigma_k}{2\gamma_{k+1}} \| {\boldsymbol y}_k \!-\! {\boldsymbol v}_{k} \|^2 \!+\! \alpha_{k} \! \Bigg( \! F\left(T_{L_k} \! ({\boldsymbol y}_k) \right) \!+\! \frac{1}{2L_k} \| r_{L_k} \! ({\boldsymbol y}_k) \|^2 \!+\! \sum_{j = 1}^{k-1} \frac{\beta_{j, k} \gamma_j}{2} \| {\boldsymbol v}_j \!-\! {\boldsymbol y}_k \|^2 \Bigg) \nonumber \\ 
			&- \frac{\alpha_k^2}{2\gamma_{k+1}} \| r_{L_k} ({\boldsymbol y}_k) \|^2 + \frac{(1-\alpha_k) \gamma_k}{2} \| {\boldsymbol x}_{\Phi_k}^* - {\boldsymbol v}_k \|^2 + \frac{(1 - \alpha_k)\alpha_k^2}{\gamma_{k+1}} \sum_{j = 1}^{k-1} \beta_{j, k} \gamma_j ({\boldsymbol v}_j - {\boldsymbol y}_k)^T r_{L_k} ({\boldsymbol y}_k) \nonumber \\ 
			& + \frac{\alpha_k^3}{\gamma_{k+1}} \sum_{j = 1}^{k-1} \beta_{j, k} \gamma_j \| {\boldsymbol v}_j - {\boldsymbol y}_k \| \| r_{L_k} ({\boldsymbol y}_k) \| + \frac{\alpha_k^2}{\gamma_{k+1} } \sum_{j = 1}^{k-1} \beta_{j, k} \gamma_j ({\boldsymbol v}_j - {\boldsymbol y}_k)^T r_{L_k} ({\boldsymbol y}_k) + \sum_{j = 1}^{k} \frac{\beta_{j, k+1} \gamma_j}{2} \| {\boldsymbol x}_{\Phi_{k+1}}^* - {\boldsymbol v}_j \|^2 \nonumber \\ 
			&+ \frac{\alpha_k (1-\alpha_k)\gamma_{k}}{\gamma_{k+1}} \Bigg( ({\boldsymbol v}_{k} - {\boldsymbol y}_k)^T r_{L_k} \left( {\boldsymbol y}_k \right) + \sum_{j = 1}^{k-1} \beta_{j, k} \gamma_j \| {\boldsymbol v}_j - {\boldsymbol y}_k \| \, \| {\boldsymbol v}_k - {\boldsymbol y}_k \| \Bigg).
		\end{align}
	\end{figure*}	
	Last, recall that we want the estimating function to be as close to the objective function as possible. Thus, we let $\theta_{k+1}^*$ equal to the upper bound obtained in \eqref{theta_kkkkk}. Letting $\phi_k^* = \theta_k^*, \forall k$ concludes the proof.  
\end{proof}

\section*{Proof of Theorem~\ref{conv_analysis_t_1}}
	\begin{proof} 
		Let us begin by setting $\Phi_0^* = F({\boldsymbol x}_0)$. Further, evaluating \eqref{phi} for $k = 0$ and ${\boldsymbol x} = {\boldsymbol x}_0$ we have: $\Phi_0({\boldsymbol x}_0) = F({\boldsymbol x}_0) + \frac{\gamma_0}{2} \| {\boldsymbol x}_0 - {\boldsymbol v}_0 \|^2$. Moreover, using the initialization ${\boldsymbol v}_0 = {\boldsymbol x}_0$ as suggested in Algorithm~\ref{FGM} we obtain $F({\boldsymbol x}_0) \leq \Phi_0^*$. Last, note that the proposed method is designed to ensure $F({\boldsymbol x}_k) \leq \Phi_k^*, \; k = 1,2, \ldots$. Applying the findings from Lemma~\ref{SFGM_lemma_1} suffices to conclude the proof. 
\end{proof}

\section*{Proof of Lemma~\ref{conv_analysis_lemma_1}}
\begin{proof}
	Let $\gamma_0 \in [0, \mu_{\hat{f}}] \cup [2\mu_{\hat{f}}, 3L_0 + \mu_{\hat{f}}]$ and apply \eqref{gamma_expr} to 
	\begin{align}
		\label{k}
		\gamma_{k+1} \! - \! \sigma_k \! &= \! (1 \! - \! \alpha_k) \gamma_{k} \! + \! \alpha_k \sigma_k \! - \! \sigma_k \! . 
	\end{align}
	Moreover, since $\lambda_0 = 1$, we can re-write \eqref{k} as
	\begin{align}
		\label{kkkk}
		\gamma_{k+1} - \sigma_k = (1 - \alpha_k) \lambda_0 \left[  \gamma_{k} - \sigma_k \right].
	\end{align}
	Substituting \eqref{gamma_expr} into \eqref{kkkk}, results in 
	\begin{align}
		\label{FGM_conv_eq_1}
		\gamma_{k+1} - \sigma_k = \lambda_{k+1}\left[\gamma_{0} - \sigma_k \right].
	\end{align}
	Next, we note that \eqref{lambda_recursive} and \eqref{alpha_k_intuition} are connected through $\alpha_k$ as follows
	\begin{align}
		\label{t}
		\alpha_k \! &= \! 1 \! - \! \frac{\lambda_{k+1}}{\lambda_k} \! = \! \sqrt{\frac{\gamma_{k+1}}{L_k}} \! =  \! \sqrt{\frac{\sigma_k}{L_k} \! +  \! \frac{\gamma_{k+1} \! -  \! \sigma_k}{L_k}}. 
	\end{align}
	Moreover, replacing \eqref{FGM_conv_eq_1} in the RHS of \eqref{t}, and making some manipulations yields
	\begin{align}
		\label{kkkkk}
		\frac{\lambda_k - \lambda_{k+1}}{\lambda_k \lambda_{k+1}} &= \frac{1}{\sqrt{\lambda_{k+1}}} \sqrt{\frac{\sigma_k}{\lambda_{k+1} L_k} + \frac{\gamma_{0} - \sigma_k}{L_k}}.
	\end{align}
	Observe that LHS of \eqref{kkkkk} can be written as $\frac{1}{\lambda_{k+1}} - \frac{1}{\lambda_k}$. Replacing the relation for the difference of squares in the LHS of \eqref{kkkkk} results in 
	\begin{align}
		\nonumber
		\left( \! \frac{1}{\sqrt{\lambda_{k+1}}} \! - \! \frac{1}{\sqrt{\lambda_{k}}} \! \right) \!\! &\left( \! \frac{1}{\sqrt{\lambda_{k+1}}} \! +  \! \frac{1}{\sqrt{\lambda_{k}}} \! \right) \!=\! \frac{1}{\sqrt{\lambda_{k+1}}} \\ 		\label{convergence_stupid}
		&\times \sqrt{ \! \frac{\sigma_k}{\lambda_{k+1} L_k} \! + \! \frac{\gamma_{0} \! - \! \sigma_k}{L_k}} \! . 
	\end{align}
	
	Observe that in Lemma \ref{SFGM_lemma_2} we define $\alpha_k \in [0,1]$. Moreover, based on \eqref{lambda_recursive} we can establish that $\lambda_k$ are non-increasing in $k$. This allows for replacing $\frac{1}{\sqrt{\lambda_{k}}}$ in the LHS of \eqref{convergence_stupid} with $\frac{1}{\sqrt{\lambda_{k+1}}}$, which would have a bigger value. So, we obtain
	\begin{align}
		\label{FGM_conv_eq_2}
		\frac{2}{\sqrt{\lambda_{k+1}}} \! \left(  \! \frac{1}{\sqrt{\lambda_{k+1}}} \! - \! \frac{1}{\sqrt{\lambda_{k}}} \! \right)  \! &\geq \! \frac{1}{\sqrt{\lambda_{k+1}}} \sqrt{\frac{\sigma_k}{\lambda_{k+1} L_k} \! + \! \frac{\gamma_{0} \! - \! \sigma_k}{L_k}}.
	\end{align}
	
	We can now observe that the convergence rate of the minimization process is dependent on the value of $\gamma_0$. We will prove convergence separately for $ \gamma_0 \in \mathcal{R}_1 = [0, \mu_{\hat{f}}[$ and $\gamma_0 \in \mathcal{R}_2 = [2\mu_{\hat{f}}, 3L_k + \mu_{\hat{f}}]$. We start with $\gamma_0 \in \mathcal{R}_1$ and introduce the following
	\begin{align}
		\label{xi_k_def}
		\xi_{k, \mathcal{R}_1} \triangleq \sqrt{\frac{L_{\max}}{\left(\sigma_k - \gamma_{0} \right) \lambda_{k}}}.
	\end{align}
	Next, we can revise \eqref{FGM_conv_eq_2} as
	\begin{align}
		\label{sfgm_conv_useless}
		\frac{2}{\sqrt{\lambda_{k+1}}} - \frac{2}{\sqrt{\lambda_{k}}}  &\geq \sqrt{\frac{\sigma_k - \gamma_{0}}{L_k}} \sqrt{\frac{\mu_{\hat{f}} L_k}{L_k \lambda_{k+1} \left(\sigma_k - \gamma_{0} \right )} + 1}.
	\end{align}
	Revising the LHS in \eqref{sfgm_conv_useless} and multiplying by $\sqrt{\frac{L_{\max}}{\sigma_k - \gamma_0}}$, yields
	\begin{align}
		\label{FGM_conv_eq_3}
		\xi_{k+1, \mathcal{R}_1} - \xi_{k, \mathcal{R}_1} &\geq \frac{1}{2}\sqrt{\frac{\sigma_k \xi_{k+1, \mathcal{R}_1}^2}{L_{\max}} + 1}.
	\end{align}
	
	Next, we prove by induction that
	\begin{align}
		\label{FGM_conv_eq_4}
		\xi_{k, \mathcal{R}_1} \geq \frac{\sqrt{2}}{4 \delta} \sqrt{\frac{L_k}{\mu_{\hat{f}} - \gamma_0}} \left[e^{(k+1) \delta} - e^{(k+1) \delta}\right],
	\end{align}
	where $\delta \triangleq \frac{1}{2} \sqrt{\frac{\sigma_k}{L_{\max}}}$. First, considering \eqref{xi_k_def} at iteration $k = 0$ and recalling that $\lambda_0 = 1$, yeids
	\begin{align}
		\label{loose_bound}
		\xi_{0, \mathcal{R}_1} &= \sqrt{\frac{L_{\max}}{(\mu_{\hat{f}} + \gamma_{-1} - \gamma_{0}) \lambda_{0}}} = \sqrt{\frac{L_{\max}}{\mu_{\hat{f}} - \gamma_{0}}}.
	\end{align}
	Embedding \eqref{L_bound} in \eqref{loose_bound}, results in 
	\begin{align}
		\xi_{0, \mathcal{R}_1} &\geq \frac{\sqrt{2}}{2} \sqrt{ \frac{L_k}{\mu_{\hat{f}} - \gamma_0}} \left[e^{\sqrt{2}/2} - e^{-\sqrt{2}/2}\right] \nonumber \\ 
		&\geq \frac{\sqrt{2}}{4 \delta} \sqrt{\frac{L_k}{\mu_{\hat{f}} - \gamma_0}} \left[e^\delta - e^{-\delta}\right] \label{tt}.
	\end{align}
	The last inequality in \eqref{tt} holds true because the RHS increases together with $\delta$, which is designed such that $\delta < \frac{\sqrt{2}}{2}$.
	
	Now suppose that \eqref{FGM_conv_eq_4} holds true at step $k$, and prove the relation for step $k+1$ by contradiction. Let $\omega(t) \triangleq \frac{\sqrt{2}}{4 \delta} \sqrt{\frac{L_k}{\mu_{\hat{f}} - \gamma_0}} \left[e^{(t+1) \delta} - e^{-(t+1) \delta}\right]$. Based on \cite[Lemma 2.2.4]{Nesterov_book} $\omega(t)$ is convex in $t$. So, we have 
	\begin{align}
		\label{FGM_conv_eq_5}
		\omega(t) \leq \xi_{k, \mathcal{R}_1} \stackrel{}{\leq} \xi_{k+1, \mathcal{R}_1} - \frac{1}{2}\sqrt{\frac{\sigma_k \xi_{k+1, \mathcal{R}_1}^2}{L_{\max}} - 1},
	\end{align}
	where the second inequality stems from \eqref{FGM_conv_eq_3}. Moreover, suppose that $\xi_{k+1, \mathcal{R}_1} < \omega(t+1)$ and substitute the relation in \eqref{FGM_conv_eq_5}. This yelds
	\begin{align}
		\omega(t) &< \omega(t+1) - \frac{1}{2}\sqrt{\frac{\sigma_k \xi_{k+1, \mathcal{R}_1}^2}{L_{\max}} - 1}.
	\end{align}
	Applying the definition for $\delta$, together with \eqref{FGM_conv_eq_4}, results in the following inequality 
	\begin{align} \label{omega(t)f}
		\omega(t) &\leq \omega(t + 1) \nonumber \\
		&- \frac{1}{2}\sqrt{4 \delta^2 \left[ \frac{\sqrt{2}}{4 \delta} \sqrt{ \frac{L_k}{\mu_{\hat{f}} - \gamma_0}} \left(e^{(t+2) \delta} - e^{-(t+2) \delta}\right) \right]^2 \! \! \! \! - \! \! 1} \nonumber \\
		&\leq \omega(t+1) - \frac{\sqrt{2}}{4} \sqrt{\!\frac{L_k}{\mu_{\hat{f}} - \gamma_0}} \left[e^{(t+2) \delta} + e^{-(t+2) \delta}\right] \\ \nonumber
		&= \omega(t+1) + \omega'(t+1)\left(t - (t+1)\right) \leq \omega(t).
	\end{align}
	The last inequality is obtained based on the supporting hyperplane theorem of convex functions. At this point, we highlight the contradiction with the earlier assumption, i.e., $\xi_{k+1, \mathcal{R}_1} < \omega(t+1)$. So, it must be true that \eqref{FGM_conv_eq_4} holds for all iterations $k = 0, 1, \ldots$. 
	
	We can now prove \eqref{FGM_conv_eq_66}. Considering \eqref{xi_k_def}, we have
	\begin{align}
		\label{jj}
		\lambda_{k} &= \frac{L_{\max}}{\xi_{k+1, \mathcal{R}_1}^2 (\sigma_k - \gamma_0)}. 
	\end{align}
	Substituting \eqref{FGM_conv_eq_4} into \eqref{jj}, yields
	\begin{align}
		\label{ttttttttttttt}
		\lambda_{k} \leq \frac{ (4 \delta)^2 L_{\max}}{2 L_k \left[e^{(k+1) \delta} - e^{(k+1) \delta}\right]^2}.
	\end{align}
	The first inequality in \eqref{FGM_conv_eq_66} is obtained by replacing the definition of $\delta$ in \eqref{ttttttttttttt}. The second inequality in \eqref{FGM_conv_eq_66} can be proved as follows. First, let us define the following abbreviation
	\begin{align}
		\mathcal{A}_k \triangleq \left(e^{\frac{k + 1}{2} \sqrt{\frac{\sigma_k}{L_k}}} \! - \! e^{-\frac{k + 1}{2} \sqrt{\frac{\sigma_k}{L_k}}}\right)^2 \!
	\end{align}
	Now, consider
	\begin{align}
		\mathcal{A}_k \label{non-strongly-cvx-lambda}
		&= \! e^{\left(k \! + \! 1\right) \sqrt{\frac{\sigma_k}{L_k}}} \! - \! e^{-\left(k+1\right) \sqrt{\frac{\sigma_k}{L_k}}} \! - \! 2. 
	\end{align}
	Applying the definition of the hyperbolic cosine function in \eqref{non-strongly-cvx-lambda}, yields
	\begin{align}
		\mathcal{A}_k = \! 2 \text{cosh}\left(\sqrt{\frac{\sigma_k}{L_k}} \left(k+1\right) \! - \! 2\right). 
	\end{align}
	Taking the Taylor expansion of $\text{cosh} (\cdot)$, yields
	\begin{align}
		\mathcal{A}_k \! &= \! -2 \! + \! 2 \! + \! 2 \frac{\sigma_k \left(k \! + \! 1\right)^2}{2L_k} \! + \! 2  \frac{\sigma_k^2 (k+1)^4}{4! {L_k}^2} + \ldots . \label{to be truncated}
	\end{align}
	Discarding the additional terms in \eqref{to be truncated} we obtain
	\begin{align}
		\label{to be substituted}
		\mathcal{A}_k \geq \frac{\sigma_k}{L_k} \left(k+1\right)^2 .
	\end{align}
	Replacing \eqref{to be substituted} in the denominator of the first inequality of \eqref{FGM_conv_eq_66} concludes the first part of the proof. The results for the case when $\gamma_0 \in \mathcal{R}_2$ can be established by following the analysis conducted for FGM in \cite[Lemma 2.2.4]{Nesterov_book}. The main update would need to be the addition of the term $\sum_{i = 1}^{k-1} \beta_{i, k} \gamma_i$ in the update for the sequence $\{\gamma_k\}_k$.
\end{proof}

\section*{Proof of Theorem~\ref{th3}}
\begin{proof}
Combining \eqref{ll} and Lemma~\ref{conv_analysis_lemma_1} for both cases of $\gamma_0 \in [0, \mu_{\hat{f}}[$ and $\gamma_0 \in [2\mu_{\hat{f}}, 3L_0 + \mu_{\hat{f}}]$ with Theorem~\ref{conv_analysis_t_1} immediately yields the convergence rates for the corresponding cases. The convergence rates in these two cases differ from each other only by a constant factor, which is $\mu_{\hat{f}} / L_k$ for $\gamma_0 \in [0, \mu_{\hat{f}}[$ and $2 \mu_{\hat{f}} / (\gamma_0 - \mu_{\hat{f}})$ for $\gamma_0 \in [2\mu_{\hat{f}}, 3L_0 + \mu_{\hat{f}}]$. It is expected that this constant facror is smaller for  $\gamma_0 \in [2\mu_{\hat{f}}, 3L_0 + \mu_{\hat{f}}]$. 
\end{proof}


\end{document}